\documentclass[a4paper,12pt,reqno]{amsart}
\usepackage[applemac]{inputenc}
\usepackage[T1]{fontenc}
\usepackage{amsfonts} 
\usepackage{latexsym} 
\usepackage{amsmath}
\usepackage{amscd}
\usepackage{ae}
\usepackage{amssymb}
\usepackage{amsthm} 
\usepackage[dvips]{graphicx}
\usepackage{longtable}
\usepackage{mathrsfs}
\usepackage[colorlinks]{hyperref}
\theoremstyle{plain}

\numberwithin{equation}{section}
\newtheorem{thm}[equation]{Theorem}    
\newtheorem{lem}[equation]{Lemma}

\newtheorem{prop}[equation]{Proposition}

\theoremstyle{definition}
\newtheorem{dfn}[equation]{Definition}

\theoremstyle{remark}
\newtheorem{rmk}[equation]{Remark}

\def\R{\mathbb{R}}

\def\Z{\mathbb{Z}}
\def\N{\mathbb{N}}

\def\viiva #1 #2 {\mathop{\Big/}_{\!\!\!{#1}}^{\>\,{#2}}}
\def\kint_#1{\mathchoice%
          {\mathop{\kern 0.2em\vrule width 0.6em height 0.69678ex depth -0.58065ex
                  \kern -0.8em \intop}\nolimits_{\kern -0.4em#1}}%
          {\mathop{\kern 0.1em\vrule width 0.5em height 0.69678ex depth -0.60387ex
                  \kern -0.6em \intop}\nolimits_{#1}}%
          {\mathop{\kern 0.1em\vrule width 0.5em height 0.69678ex depth -0.60387ex
                  \kern -0.6em \intop}\nolimits_{#1}}%
          {\mathop{\kern 0.1em\vrule width 0.5em height 0.69678ex depth -0.60387ex
                  \kern -0.6em \intop}\nolimits_{#1}}}

\begin{document}
\author{Lizaveta Ihnatsyeva AND Antti V. V\"{a}h\"{a}kangas}

\title[Characterization of traces of smooth functions]{Characterization of traces 
of smooth functions on Ahlfors regular sets}

\keywords{Triebel--Lizorkin space, Besov space, 
Ahlfors $d$-regular set, trace theorem, local polynomial approximation}
\subjclass[2000]{46E35}

\begin{abstract}
We extend the results of P. Shvartsman on characterizing the traces
of Besov and Triebel--Lizorkin spaces on Ahlfors $n$-regular sets to the case of
$d$-regular sets, $n-1<d<n$. 
The characterizations of trace spaces are given in terms of local polynomial approximations.
\end{abstract}

\maketitle

\section{Introduction}

The problem of characterizing the traces of commonly used spaces of smooth functions on different subsets of $\R^n$ has been 
extensively studied by several authors, see e.g. the monographs \cite{BesovIlinNikolski,JonssonWallin1984,MazyaPoborchi,Triebel2,Triebel3} and the references therein. In this paper, we focus on a certain description of the traces of Besov spaces and Triebel--Lizorkin spaces on Ahlfors $d$-regular sets, where $n-1<d<n$. 
In particular, since Triebel--Lizorkin spaces include Sobolev spaces
$W^{k,p}(\R^n)$, $k \in\N$, $1<p<\infty$, we cover trace theorems for these spaces as well.

The first approach to the notion of Besov spaces on $d$-sets, $0<d\le n$, was developed in 
a number of papers by A. Jonsson and H. Wallin, see \cite{JonssonWallin1984} for a unified treatment of this material. Their definition of Besov spaces is based on a jet technique and 
it gives a characterization for the traces of the classical Besov spaces on $d$-regular subsets of $\R^n$, in the sense that the trace operator is a bounded linear surjection  with a 
bounded linear right inverse (the extension operator). When $d<n$ the Besov spaces introduced by  A. Jonsson and H. Wallin also describe the trace spaces for Bessel potential spaces \cite{JonssonWallin1984} 
and for their natural extensions, Triebel--Lizorkin spaces, see 
the papers of G. Mamedov \cite{Mamedov} and K. Saka \cite{SakaI}, \cite{SakaII}.

Our work is motivated by
the results of P. Shvartsman, who has earlier given an intrinsic characterization of the traces of Besov spaces and 
Triebel--Lizorkin spaces on Ahlfors $n$-regular subsets of $\R^n$ via local polynomial approximations \cite{Shvartsman}. 
In this paper, we extend the
results in \cite{Shvartsman} to the case of Ahlfors $d$-regular sets with $n-1<d<n$. This possibility was briefly mentioned by Yu. Brudnyi in the survey paper \cite{Brudnyi}.

The definition of Besov spaces on $d$-sets and the related methods, based on jet technique, are rather technical in general. 
Using a local polynomial approximation approach we define the Besov spaces on a $d$-set,  $n-1<d<n$, and we show that these spaces can be equipped
with a wide family of equivalent norms, see Section \ref{besov_sec}
and Theorem \ref{beseq}.
The 
definition we use has a simpler formulation and 
it allows to simplify proofs of the trace theorems.

In Section \ref{restr_sec} we prove the restriction parts of the trace theorems for Besov and Triebel--Lizorkin spaces, i.e.
boundedness of the pointwise trace operator whose range coincides with a Besov space on a $d$-set in both cases,
see Theorem \ref{restrB} and Theorem \ref{ctrace}.
Let us note that in case of 
Triebel--Lizorkin spaces the crucial argument in the proof is that $d$-sets with $d<n$ are porous (porosity is also known as the ball condition), 
a fact which fails for $n$-sets.

The extension part of the trace theorems is treated in sections \ref{ext_sec} and \ref{ext_th_sec}, Theorem \ref{bes_est} and Theorem \ref{thf} being our extension theorems. In
Section \ref{ext_sec} we construct linear
extension operators by generalizing the construction in \cite{Shvartsman}. We rely on a Remez type inequality from \cite{BrudnyiBrudnyi} which leads
us to the assumption $n-1<d$, a lower bound for the dimension of a set. The proof of the boundedness of the extension operator again uses the fact that $d$-sets with $d<n$ are porous,
more precisely, we need estimates established in Section \ref{prel_sec}.
The significance of porosity in the context of trace theorems was noticed already by M. Frazier and B. Jawerth \cite[Theorem 13.7]{fj}, see also 
\cite[Theorem 9.21]{Triebel3} in the monograph of H. Triebel.

\section{Preliminaries}\label{prel_sec}

\subsection{Notation}
Let $x,y\in\R^n$, then $\|x-y\|_\infty=\max_{i=1,\ldots,n} |x_i-y_i|$ denotes the distance
between $x$ and $y$ in the supremum metric. The distance
between a point $a\in\R^n$ and a set $\emptyset\not=X\subset\R^n$ in the
supremum metric is denoted by
$\mathrm{dist}(a,X)$, 
the diameter of $X$ in the supremum metric is $\mathrm{diam}\,X$.

By $Q=Q(x_Q,r_Q)$ we denote a closed cube in $\mathbb{R}^n$ centered at $x_Q\in \R^n$ with the
side length $\ell(Q)=2r_Q$ and sides parallel to the coordinate axes. We write $tQ$, $t>0$, for the 
cube which is centered 
at $x_Q$ and whose
side length is $t\ell(Q)$. We denote by
$\mathcal{D}$ the family of closed dyadic cubes in $\R^n$; $\mathcal{D}_j$ stands for
the family of those dyadic cubes whose
side length is $2^{-j}$, $j\in\Z$.

We denote by $|x-y|_2$ the Euclidean distance between  $x$ and $y$ in $\R^n$. The Euclidean diameter of a set $X$ in $\R^n$
is 
$\mathrm{diam}_2 X$.
We write $B(x,r)$ for an open ball in $\R^n$ with the center $x$ and the
radius $r>0$. If $x\in\R^n$ and $c>0$, then
$cB(x,r)$ stands for the ball $B(x,cr)$.

The Lebesgue measure of a measurable set $E$ in $\R^n$ is denoted by $|E|$.
By $\mathcal{H}^d$ we denote the $d$-dimensional Hausdorff measure in $\R^n$.
For an $\mathcal{H}^d$-measurable set $E$ with a finite and positive measure, we denote
$$\kint_Ef\,d\mathcal{H}^d=\frac{1}{\mathcal{H}^d(E)}\int_{E}f\,d\mathcal{H}^d.$$
We write $\chi_E$ for the characteristic function of a set $E$.
The symbol $c$ is used for various positive constants; they may change even on the same line. The dependence on parameters is expressed, for example, by $c=c_{S,n,k}$.

\subsection{Function spaces in $\R^n$}\label{SecFunctionSpacesRn}
First recall the definitions and basic properties of Sobolev spaces, Besov spaces and Triebel--Lizorkin spaces on $\R^n$. 

Let $L^p(\R^n)$ denote the Lebesgue space of $p$-integrable functions in $\R^n$. We write \[\|f\|_p=\|f\|_{L^p(\R^n)}:=\bigg(\int_{\R^n} |f(x)|^p dx\bigg)^{1/p}.\]
For a positive integer $k$ and $1\leq p<\infty$, 
the Sobolev space $W^{k,p}(\R^n)$ consists of all functions $f\in L^p(\R^n)$ 
having  distributional derivatives $\partial^\alpha f$,
$|\alpha|\leq k$, in $L^p(\R^n)$. The Sobolev space is equipped with a norm
$\|f\|_{W^{k,p}(\R^n)} := \sum_{|\alpha|\le k}\|\partial^\alpha f\|_{p}$.

There are several equivalent characterizations for the Sobolev spaces and their natural extensions, Besov spaces and Triebel--Lizorkin spaces (see e.g. \cite{AdamsHedberg}, \cite{BesovIlinNikolski}, \cite{Triebel89} and \cite{Triebel2} for the general theory). In this
paper, we mostly use definitions based on the local polynomial approximation approach.

Let $f\in L^u_{\mathrm{loc}}(\R^n)$, $1\le u \le \infty$, and $k\in\N_0$.
The {\em normalized local best approximation} of $f$ on a cube
$Q$ in $\R^n$ is defined by
\begin{equation} 
\mathcal{E}_k(f,Q)_{L^u(\R^n)}:=\inf_{P\in \mathcal{P}_{k-1}}\bigg(\frac{1}{|Q|}\int_{Q}|f(x)-P(x)|^u\;dx\bigg)^{1/u}.
\end{equation}
Here and below  $\mathcal{P}_{k}$, $k\geq 0$, denotes the space of polynomials on $\mathbb{R}^n$ of degree at most $k$; we also write $\mathcal{P}_{-1}:=\{0\}$ for convenience.

For $\alpha>0$ and $1\le p,q\le\infty$, the Besov space $B_{pq}^{\alpha}(\R^n)$ can be defined in the following way, see e.g. \cite{Brudnyi74}, \cite{Triebel89}, \cite{Triebel2}.
Let $k$ be an integer such that $\alpha<k$ and $1\le u\le p$, then $B_{pq}^{\alpha}(\R^n)$ consists of functions $f\in L^p(\R^n)$ such that
\begin{equation}\label{DefofBesovSpace}
\int_0^1\bigg(
\frac{\Vert\mathcal{E}_k(f,Q(\cdot,t))_{L^u(\R^n)}\Vert_{L^p(\R^n)}}{t^{\alpha}}\bigg)^{q}\;\frac{dt}{t}<\infty,\qquad 
\text{ if }q<\infty,
\end{equation}
and $\sup_{0<t\le 1} t^{-\alpha}\Vert\mathcal{E}_k(f,Q(\cdot,t))_{L^u(\R^n)}\Vert_{L^p(\R^n)}<\infty$ if $q=\infty$. The Besov norms
$$
\|f\|_{B^\alpha_{pq}(\R^n)}:=\|f\|_{L^p(\R^n)}+\bigg(\int_0^1\bigg(
\frac{\Vert\mathcal{E}_k(f,Q(\cdot,t))_{L^u(\R^n)}\Vert_{L^p(\R^n)}}{t^{\alpha}}\bigg)^{q}\;\frac{dt}{t}\bigg)^{1/q}
$$
(modification
if $q=\infty$) are equivalent if $1\le u\le p$ and $\alpha<k$.

The following definition of Triebel--Lizorkin spaces can be found e.g.
in \cite{Dorronsoro} or \cite{Triebel89}.
Let $\alpha>0$, $1<p<\infty$, $1<q\le\infty$ and $k$ be an integer such that $\alpha<k$. For $f\in L^1_{\mathrm{loc}}(\R^n)$ and $x\in \R^n$, we
set
\[
g(x):=\bigg(\int_0^1\bigg(\frac{\mathcal{E}_k(f,Q(x,t))_{L^1(\R^n)}}{t^{\alpha}}\bigg)^{q}\;\frac{dt}{t}\bigg)^{1/q} ,\qquad \text{if }q<\infty,
\]
and $g(x):=\sup\{t^{-\alpha}\mathcal{E}_k(f,Q(x,t))_{L^1(\R^n)}:\, 0<t\le 1\}$ if $q=\infty$. 
The function $f$ belongs to
Triebel--Lizorkin space $F^\alpha_{pq}(\R^n)$ if $f$ and $g$ are both in $L^p(\R^n)$. The Triebel--Lizorkin norms
$$
\|f\|_{F^\alpha_{pq}(\R^n)}:=\|f\|_{L^p(\R^n)}+\|g\|_{L^p(\R^n)}
$$
are equivalent if $\alpha<k$.

Triebel--Lizorkin spaces include Sobolev spaces as a special case:
 $F_{p2}^{k}(\R^n)$ coincides with $W^{k,p}(\R^n)$ for any $k\in\N$ and $1<p<\infty$, see e.g. \cite{AdamsHedberg} for a proof.

\subsection{Ahlfors regular sets}
Recall also the definition of an Ahlfors $d$-regular set,
or, briefly,  $d$-set.

\begin{dfn}
Let $0<d\le n$.
A closed set $S\subset \R^n$ is a {\em $d$-set}  if there are constants $c_1,c_2>0$ such that
\begin{equation}\label{sset}
c_1r^d \le \mathcal{H}^d(Q(w,r)\cap S)\le c_2r^d
\end{equation}
for every $w\in S$ and $0<r\le 1$.
\end{dfn}

Note that if $S$ is a $d$-set then the estimate
\begin{equation}\label{sset_rem}
c^{-1}r^d\le \mathcal{H}^d(Q(w,r)\cap S)\le cr^d 
\end{equation}
holds true for every $w\in S$ and $0<r\le R$, where $R$ is any fixed positive number and the constant $c\ge 1$ depends on
parameters $R,c_1,c_2,d,n$.

\begin{rmk}
Sometimes the definition of a $d$-set is formulated in terms of a Borel measure $\mu$ in $\R^n$, satisfying
the conditions $\mathrm{supp}\,\mu =S$ and
\eqref{sset} with $\mathcal{H}^d$ replaced by $\mu$. Such a measure $\mu$ can be identified, up to constants $c_1$ and $c_2$, with the measure
$\mathcal{H}^d|_S$; a short proof can be found in \cite[Theorem 3.4]{Triebelfs}.
This equivalence of measures also
implies that \[\mathcal{H}^d(\bar S\setminus \mathrm{int}\,S)=0\] (see \cite[pp. 205--206]{JonssonWallin1984}) and, therefore, we do not lose generality when we consider only closed sets. 
\end{rmk}

\subsection{Porous sets}\label{application}
In this paper we widely use the fact that $d$-sets with $d<n$ are porous (porosity is also known as the ball condition). Here we prove a 
useful norm estimate for porous sets, Theorem \ref{reverse}.

\begin{dfn}\label{porous}
A  set $S\subset\R^n$ is 
{\em porous} (or {\em $\kappa$-porous}) if for some $\kappa\ge 1$
the following statement is true:
For every cube $Q(x,r)$ with $x\in\R^n$ and $0<r\le 1$ there is $y\in Q(x,r)$ such that
$Q(y,r/\kappa)\cap S=\emptyset$.
\end{dfn}

\begin{rmk}\label{porous_lem}
The observation that $d$-sets with $d<n$ are porous was already done in \cite{JonssonII}. See also Proposition 9.18 in \cite{Triebel3} which gives this fact as a special case.
Let us also mention that a set $S$ in $\R^n$
is porous if, and only if, its Assouad dimension strictly less 
than $n$, \cite{Luukkainen}. 
\end{rmk}

For a  $S$ in $\R^n$ and positive constant $\gamma>0$ we denote
\begin{equation}\label{c_definition}
\mathcal{C}_{S,\gamma} = \{Q\in\mathcal{D}\,:\,\gamma^{-1}\mathrm{dist}(x_Q,S)\le \ell(Q)\le 1\}.
\end{equation}
This is the family of dyadic cubes that are relatively close to the set.

\begin{thm}\label{reverse}
Suppose that $S\subset\R^n$ is porous.
Let $p,q\in (1,\infty)$ and $\{a_Q\}_{Q\in\mathcal{C}_{S,\gamma}}$ be a sequence of non-negative scalars. 
Then
\begin{equation}\label{rev}
\bigg\|\sum_{Q\in\mathcal{C}_{S,\gamma}} \chi_Q a_Q\bigg\|_p\le c\bigg\|\bigg(\sum_{Q\in\mathcal{C}_{S,\gamma}}
(\chi_{Q} a_Q)^q\bigg)^{1/q}\bigg\|_p.
\end{equation}
Here the constant $c$ depends on $n$, $p$, $\gamma$ and the set $S$.
\end{thm}

The proof of this theorem
is a consequence of maximal-function techniques that can be found in \cite{Bojarski} and \cite{iwaniec}. 
First we consider some auxiliary statements which are needed for the proof.
Let $S\subset \R^n$ be a $\kappa$-porous set, $\gamma>0$ and $\mathcal{C}_{S,\gamma}$ be 
as in \eqref{c_definition}. Define
\begin{equation}\label{sigdef}
\sigma=(1+\gamma)16\kappa>1.
\end{equation}

The following lemma assigns to every $Q\in\mathcal{C}_{S,\gamma}$
a dyadic cube  $r(Q)\subset Q$ that is away from $S$
and 
 the measures of $Q$ and $r(Q)$
are comparable.

\begin{lem}\label{kootut}
Let
$Q\in\mathcal{C}_{S,\gamma}$. Then there is  $r(Q)\in\mathcal{D}$ such that
$\ell(Q)\le \sigma\ell(r(Q))$ and
\[
r(Q)\subset \mathrm{int}(Q)\cap \{x\in\R^n\,:\,\ell(Q)/\sigma\le \mathrm{dist}(x,S)\le \sigma\ell(Q)\}.
\]
\end{lem}

\begin{proof}
Let $Q=Q(x_Q,r_Q)\in\mathcal{C}_{S,\gamma}$. Then $r_Q\le 1/2$ and by
the definition of porosity there is a point $y\in Q(x_Q,r_Q/2)$ such that
$Q(y,r_Q/2\kappa)\cap S=\emptyset$.
Let $r(Q)$ be the dyadic cube
for which $y\in r(Q)$ and \[r_Q/8\kappa< \ell(r(Q))\le r_Q/4\kappa.\]
It is straightforward to verify 
that $r(Q)$ satisfies the required properties.
\end{proof}

The following proposition can be found 
in \cite[Lemma 4.2]{Bojarski} or \cite[Lemma 4]{iwaniec}, but we recall its proof
for the convenience.

\begin{prop}\label{separate}
Let $p\in (1,\infty)$ and $\{a_Q\}_{Q\in\mathcal{C}_{S,\gamma}}$ be a sequence of non-negative scalars. 
Then
\begin{equation}\label{ps}
\bigg\| \sum_{Q\in\mathcal{C}_{S,\gamma}} \chi_Qa_Q\bigg\|_p\le c\bigg\|\sum_{Q\in\mathcal{C}_{S,\gamma}}
\chi_{r(Q)}a_Q\bigg\|_p.
\end{equation}
Here $c$ depends on $n,p$ and $\sigma$.
\end{prop}

\begin{proof}
By duality it suffices to verify that
\begin{equation}\label{dual}
I:=\int_{\R^n} \Big|\sum_{Q\in\mathcal{C}_{S,\gamma}}a_Q \chi_Q(y) \psi(y)\Big|dy\le c\bigg\|
\sum_{Q\in\mathcal{C}_{S,\gamma}} a_Q \chi_{r(Q)}\bigg\|_{p}
\end{equation}
for every $\psi\in C^\infty_0(\R^n)$ with $\|\psi\|_{p'}\le 1$ ($1/p+1/p'=1$). Let
$\psi$ be such a test function and
$Q\in \mathcal{C}_{S,\gamma}$. Then
\begin{equation}\label{laste}
\int_Q |\psi(y)|dy\le |Q|M\psi(x),\qquad x\in Q.
\end{equation}
Here $M$ stands
for the non-centered Hardy--Littlewood maximal operator
over cubes.
By averaging the inequality \eqref{laste} over 
cube $r(Q)\subset Q$, and using Lemma \ref{kootut}, we get
\[
\int_Q |\psi(y)|dy \le \frac{|Q|}{|r(Q)|} \int_{r(Q)} M\psi(x)dx\le \sigma^n\int_{r(Q)}M\psi(x)dx,\quad Q\in\mathcal{C}_{S,\gamma}.
\]
By the triangle-inequality and the previous estimates
\begin{align*}
I&\le \sum_{Q\in\mathcal{C}_{S,\gamma}} a_Q\int_Q |\psi(y)|dy\\&\le \sigma^n\sum_{Q\in \mathcal{C}_{S,\gamma}} a_Q \int_{r(Q)} M\psi(x)dx=\sigma^n\int_{\R^n} \sum_{Q\in\mathcal{C}_{S,\gamma}} a_Q \chi_{r(Q)}(x) M\psi(x)dx.
\end{align*}
In order to obtain the required estimate \eqref{dual} it suffices to apply the H\"older's inequality to
the right hand side and use the boundedness of $M$ in $L^{p'}$.
\end{proof}

Observe that the family of cubes
$\{r(Q)\,:\,Q\in\mathcal{C}_{S,\gamma}\}$
is locally finite.
This follows from Lemma \ref{kootut} and the fact
that the dyadic cubes with equal side length have disjoint interiors.
Next we
partition the family $\{r(Q)\}$ into a finite number of disjoint subfamilies in the following way. 

Let $r_0\in\N$ be such that $\sigma^2< 2^{r_0}$. For every $r\in\{0,1,\ldots,r_0-1\}$, we denote
\[
\mathcal{C}_{S,\gamma}^r=\mathcal{C}_{S,\gamma}\cap \bigcup_{\substack{j\in\Z\\j\equiv r\,\mathrm{mod}\,r_0}} \mathcal{D}_j.
\]

\begin{lem}\label{big}
Let $r\in \{0,1,\ldots,r_0-1\}$ and $Q,R\in\mathcal{C}_{S,\gamma}^r$ be such that
$R\not=Q$. Then  $r(Q)\cap r(R)=\emptyset$.
\end{lem}

\begin{proof}
Let $Q$, $R\in\mathcal{C}_{S,\gamma}^r$. If $\ell(Q)=\ell(R)$, then by Lemma \ref{kootut} and properties of dyadic cubes
$r(Q)\cap r(R)\subset \mathrm{int}(Q)\cap \mathrm{int}(R)=\emptyset$. 

On the other hand, if $\ell(Q)\not=\ell(R)$ then, say, 
$2^{-i}=\ell(Q)<\ell(R)=2^{-j}$ and $j<i$. By Lemma \ref{kootut}
\[r(Q)\subset  \{x\in\R^n\,:\,\mathrm{dist}(x,S)\le \sigma 2^{-i}\}\]
and
\[r(R)\subset \{x\in\R^n\,:\,2^{-j}/\sigma\le \mathrm{dist}(x,S)\}.\]
Hence, $r(Q)\cap r(R)=\emptyset$ if $\sigma 2^{-i}<2^{-j}/\sigma$, that is,
if $\sigma^2<2^{i-j}$. This is true since $\sigma^2< 2^{r_0}\le 2^{i-j}$
(notice that $i,j\equiv r\,\mathrm{mod}\,r_0$, so $i-j\ge r_0$). 
\end{proof}

Now we are ready for the proof of inequality \eqref{rev}.

\begin{proof}[Proof of Theorem \ref{reverse}]
By Proposition \ref{separate} and triangle-inequality
\begin{equation}\label{apus}
\begin{split}
\bigg\| \sum_{Q\in\mathcal{C}_{S,\gamma}} \chi_Q a_Q\bigg\|_p&\le 
c\bigg\| \sum_{Q\in\mathcal{C}_{S,\gamma}} \chi_{r(Q)} a_Q\bigg\|_p\\&\le c\sum_{r=0}^{r_0-1}\bigg\|\sum_{Q\in\mathcal{C}_{S,\gamma}^r}
\chi_{r(Q)} a_Q\bigg\|_p.
\end{split}
\end{equation}
Fix $r\in \{0,1,\ldots,r_0-1\}$. Lemma \ref{big} implies that
the cubes $\{r(Q)\}_{Q\in\mathcal{C}_{S,\gamma}^r}$ are disjoint. Therefore,
\begin{align*}
\sum_{Q\in\mathcal{C}_{S,\gamma}^r}
\chi_{r(Q)} a_Q=\bigg(\sum_{Q\in\mathcal{C}_{S,\gamma}^r} (\chi_{ r(Q)} a_Q)^q\bigg)^{1/q}.
\end{align*}
By Lemma \ref{kootut}, we have $\chi_{r(Q)}\le \chi_Q$ if $Q\in\mathcal{C}_{S,\gamma}^r$. Hence,
\begin{align*}
\sum_{Q\in\mathcal{C}_{S,\gamma}^r}
\chi_{r(Q)} a_Q&\le \bigg(\sum_{Q\in\mathcal{C}_{S,\gamma}^r} (\chi_{Q} a_Q)^q\bigg)^{1/q}
\\&\le \bigg(\sum_{Q\in\mathcal{C}_{S,\gamma}} (\chi_{Q} a_Q)^q\bigg)^{1/q}.
\end{align*}
Taking $L^p$-norms and applying \eqref{apus} we get the
required estimate \eqref{rev}.
\end{proof}

\subsection{Hardy-type inequalities}
Recall also some inequalities of Hardy type which we often use in this paper.

\begin{lem}\label{LemmaHardyInequality}
Let $0<p<\infty$ and $a_j\geq 0$, $j=0,1,\dots$. Then
\begin{equation}\label{leindler}
\sum_{j=0}^\infty 2^{\sigma j}\bigg(\sum_{i=0}^j a_i\bigg)^p\le c\sum_{j=0}^\infty 2^{\sigma j}a_j^p \,\,\,\,\,\,\text{for}\,\,\sigma<0
\end{equation}
and
\begin{equation}\label{Hardy2}
\sum_{j=0}^\infty 2^{\sigma j}\bigg(\sum_{i=j}^\infty a_i\bigg)^p\le c\sum_{j=0}^\infty 2^{\sigma j}a_j^p\,\,\,\,\,\,\text{for}\,\,\sigma>0.
\end{equation}
\end{lem}

For $p>1$ the inequalities follow from Leindler's results in \cite{Leindler}. The case $0<p\le 1$
 follows by applying the inequalities
$$\bigg(\sum_{i=0}^j a_i\bigg)^p \le \sum_{i=0}^j a_i^p,\qquad \bigg(\sum_{i=j}^\infty a_i\bigg)^p \le \sum_{i=j}^\infty a_i^p.$$

\section{Besov spaces on $d$-sets}\label{besov_sec}
In the section we always assume that $S$ is a $d$-set with $n-1<d\le n$. We define Besov spaces on such $d$-sets
and establish their independence
on certain parameters which arise
in the definition.

We write $L^p(S,\mathcal{H}^d)$ or $L^p(S)$ for the space of $p$-integrable functions on a 
$d$-set $S$ with respect to the natural Hausdorff measure $\mathcal{H}^d|_S$.

Let $f\in L^u_{\mathrm{loc}}(S,\mathcal{H}^d)$ and $Q$ be a cube in $\mathbb{R}^n$ 
centered at $S$. Then by $\mathcal{E}_k(f,Q)_{L^u(S)}$, $u\in [1,\infty)$, we denote the normalized
local best approximation of $f$ on $Q$ with respect to the $L^u(S)$-norm, i.e.
$$
\mathcal{E}_k(f,Q)_{L^u(S)}:=\inf_{P\in \mathcal{P}_{k-1}}\bigg(\frac{1}{\mathcal{H}^d(Q\cap S)}\int_{Q\cap S}|f(x)-P(x)|^u\;d\mathcal{H}^d(x)\bigg)^{1/u}.
$$
Let $Q_1\subset Q_2$ be two cubes 
centered in $S$ such that $r_{Q_2}\le R$
for some $R>0$. Then
the definition of local best approximation and
\eqref{sset_rem} yield
\begin{equation}\label{eqMonotonyOfLocalApproxSset}
\mathcal{E}_k(f,Q_1)_{L^u(S)}
\le C\bigg(\frac{r_{Q_2}}{r_{Q_1}}\bigg)^{d/u}
\mathcal{E}_k(f,Q_2)_{L^u(S)}.
\end{equation}
Here the constant $C>0$ depends on $S$ and $R$.

Following the definition
of Besov spaces via local polynomial approximations in the case of Euclidean space (see Paragraph \ref{SecFunctionSpacesRn} and the references therein) and in the case of $n$-set, see \cite{Brudnyi76,Shvartsman},  we define Besov spaces on $d$-sets as follows.

\begin{dfn}\label{besov}
Let $n-1<d\le n$
and $S$ be a $d$-set in $\R^n$.
The Besov space $B^\alpha_{pq}(S)$, $\alpha>0$, $1\le p,q\le\infty$, is the set of those functions 
$f\in L^p(S)$ for which the norm $\|f\|_{B_{pq}^{\alpha}(S)}$ is finite. Here
\begin{equation}\label{DefinitionBesovNormIntegral}
\|f\|_{B_{pq}^{\alpha}(S)}:= \|f\|_{L^p(S)}+\bigg(\int_0^1\bigg(
\frac{\Vert\mathcal{E}_k(f,Q(\cdot,\tau))_{L^u(S)}\Vert_{L^p(S)}}{\tau^{\alpha}}\bigg)^{q}\;\frac{d\tau}{\tau}\bigg)^{\frac{1}{q}}, 
\end{equation}
where $1\le u\le p$ and $k$ is an integer such that $\alpha<k$.
\end{dfn}

\begin{rmk}\label{diskr}
By the monotonicity property \eqref{eqMonotonyOfLocalApproxSset} of local approximations it easily follows that in \eqref{DefinitionBesovNormIntegral} we can replace the integral by sum, i.e. if $m\in \Z$ then
the quantity on the right hand side of \eqref{DefinitionBesovNormIntegral} is equivalent to 
\begin{equation}
\Vert f\Vert_{L^p(S)}+\bigg(\sum_{j=m}^\infty \big(2^{j \alpha}\|\mathcal{E}_{k}(f,Q(\cdot,2^{-j}))_{L^u(S)}\|_{L^p(S)}\big)^q\bigg)^{1/q}.
\end{equation}
The equivalence constant is independent
of $f$, but it can depend on $m$.
\end{rmk}

The definition of Besov spaces on $d$-sets, $n-1<d\le n$,
introduced by A. Jonsson and H. Wallin in \cite[p.135]{JonssonWallin1984}, is  equivalent to Definition \ref{besov} with parameters $k=[\alpha]+1$ and $u=p$. It was shown explicitly in \cite{IhnatsyevaKorte}. The following
theorem states that, in fact, the space given by Definition \ref{besov}
is independent
of the admissible parameters $k$ and $u$.

\begin{thm}\label{beseq}
The Besov space $B_{pq}^{\alpha}(S)$ with $\alpha>0$ and $1\le p,q< \infty$ does
not depend on the parameters $k>\alpha$ 
and $1\le u\le p$. Furthermore,
the norms corresponding to any two pairs of such parameters are equivalent.
\end{thm}

Regarding the proof of this theorem, the independence of the norms on $k$ will be established at the end of this section, see Proposition \ref{noreq}. The independence on
the parameter $u$ is a consequence
of the corresponding result in $\R^n$ and our estimates
involved in the proof of the trace theorem; see Proposition \ref{equirem} for the details.

Theorem \ref{beseq} allows us to choose the most convenient 
values of the parameters $k>\alpha$ and $u\in [1,p]$ while formulating our main results.
We would like to emphasize that this theorem
is not applied in the subsequent proofs, so there is no circular reasoning.

\subsection{Remez type theorem}
The main reason for the restriction $d>n-1$, a lower bound on the dimension of a set $S$, is that we widely use the following Remez type theorem, for
a proof, see \cite{BrudnyiBrudnyi}. Note that the statement of Theorem \ref{remez} fails when $d\le n-1$. 

\begin{thm}\label{remez}
Assume that $U\subset\R^n$ is a bounded open set and $\omega\subset U$
is such that for some $a>0$ and $d$, $n-1<d\le n$, it satisfies
\[
\mathcal{H}^d(B(x,r)\cap \omega)\le ar^d,\quad\,\,\,x\in\omega,\,\,\,
0<r\le \mathrm{diam}_2\omega.
\]
Assume also that
$$\lambda:=\frac{\{\mathcal{H}^d(\omega)\}^{n/d}}{|U|}>0.$$
Then there is a positive constant $C$ such that, for every polynomial $p$ of degree 
at most $k$,
\begin{equation}\label{RemezTypeIneq}
\bigg(\frac{1}{|U|}\int_{U}|p|^r\;dx\bigg)^{1/r}\le C
\bigg(\frac{1}{\mathcal{H}^d(\omega)}\int_{\omega}|p|^u\;d\mathcal{H}^d\bigg)^{1/u}.
\end{equation}
Here $0<u,r\le\infty$ and the constant $C$
depends on $U,\,n,\,u,\,r,\,d,\,k,\,a,\,\lambda$, and it is increasing 
in $1/\lambda$.
\end{thm}

We use the following slight modification of Theorem \ref{remez}.
\begin{prop}\label{qremez}
Let $S$ be a $d$-set, $n-1 < d\le n$. Suppose that $Q=Q(x_Q,r_Q)$ and $Q'=Q(x_{Q'},r_{Q'})$ 
are cubes  in $\R^n$ such that 
$x_{Q'}\in S$, $Q'\subset Q$ and \[0<r_Q\le R\,r_{Q'}\le R^2\] for some $R>0$. Then,
for every polynomial $p$ of degree at most $k$, we have
\[
\bigg(\frac{1}{|Q|}\int_{Q}|p|^r\;dx\bigg)^{1/r}\le C
\bigg(\frac{1}{\mathcal{H}^d(Q'\cap S)}\int_{Q'\cap S}|p|^u\;d\mathcal{H}^d\bigg)^{1/u},
\]
where $1\le u,r\le\infty$ and the constant $C$
depends on $S,\,R,\,n,\,u,\,r,\,k$.
\end{prop}

\begin{proof}
Replacing $Q$ by $2Q$ if necessary we can assume that $Q'\subset{\rm int}\,Q$. 
Denote
\begin{align*}
U&:=\mathrm{int}\,Q(0,1)
=\bigg\{\frac{y-x_Q}{r_Q}\,:\,y\in\mathrm{int}\,Q\bigg\},\\
\omega&:=
\bigg\{\frac{y-x_Q}{r_Q}\,:\,y\in Q'\cap S\bigg\}\subset U. \hspace{2.9cm}
\end{align*}
It is easy to see that $U$ and $\omega$ satisfy the assumptions of Theorem \ref{remez}. Indeed,
if $x\in \omega$ and $0<r\le \mathrm{diam}_2\omega=\mathrm{diam}_2(Q'\cap S)/r_Q$, then 
by the translation invariance and scaling property
$\mathcal{H}^d(\sigma E)=\sigma^d\mathcal{H}^d(E)$
of the Hausdorff measure  
\begin{align*}
\mathcal{H}^d(B(x,r)\cap \omega)
=r_Q^{-d}\mathcal{H}^d(B(r_Q x+x_Q,r_Q r)\cap (Q'\cap S)).
\end{align*}
Since $r_Qx+x_Q\in Q'\cap S$ as $x\in\omega$ and
$0<r_Q r\le \mathrm{diam}_2(Q'\cap S)\le 2R\sqrt n$, by
\eqref{sset_rem}
\[
\mathcal{H}^d(B(r_Q x+x_Q,r_Q r)\cap (Q'\cap S))\le \mathcal{H}^d(Q(r_Q x+x_Q,r_Q r)\cap S)\le c(r_Qr)^d.
\]
Hence, $\mathcal{H}^d(B(x,r)\cap \omega)\le c r^d$. Consider now 
\[
\lambda:=\frac{\{\mathcal{H}^d(\omega)\}^{n/d}}{|U|}=
2^{-n}\bigg\{\bigg(\frac{1}{r_Q}\bigg)^d \mathcal{H}^d(Q'\cap S)\bigg\}^{n/d}\ge c>0.
\]
Here the lower bound for $\lambda$ follows from \eqref{sset_rem} and the assumption that $r_Q\le R r_{Q'}$.

By a change of variables and Theorem \ref{remez}, we obtain
\begin{align*}
\bigg(\frac{1}{|Q|}\int_{Q}|p(x)|^r\;dx\bigg)^{1/r}
&=
\bigg(2^{-n}\int_{U}|p(r_Q x+x_Q)|^r\;dx\bigg)^{1/r}\\&\le C
\bigg(\frac{1}{\mathcal{H}^d(\omega)}\int_{\omega}|p(r_Q x+x_Q)|^u\;d\mathcal{H}^d(x)\bigg)^{1/u}
\\&=C\bigg(\frac{1}{\mathcal{H}^d(Q'\cap S)}\int_{Q'\cap S}|p|^u\;d\mathcal{H}^d\bigg)^{1/u},
\end{align*}
where the constant $C$ is independent of a cube $Q$.
\end{proof}

\begin{rmk}
Let $Q$ be a cube in $\R^n$ which is centered in $S$ and satisfies $0<\ell(Q)\le R$.
Proposition \ref{qremez} implies
the following reverse
H\"older inequality for polynomials $p\in\mathcal{P}_k$,
\begin{equation}\label{RevHolderForS-set}
\bigg(\frac{1}{\mathcal{H}^d(Q\cap S)}\int_{Q\cap S}|p|^r\;d\mathcal{H}^d\bigg)^{1/r}\le C\bigg(\frac{1}{\mathcal{H}^d(Q\cap S)}\int_{Q\cap S}|p|^u\;d\mathcal{H}^d\bigg)^{1/u}.
\end{equation} 
Here $1\le u,r\le\infty$ and the constant $C$ depends on $S,\,R,\,n,\,u$  and $k$.
\end{rmk}

\subsection{Projections to polynomials}\label{projections_sec}
Here we construct certain polynomial projection operators; the construction is an extension of the case $d=n$ in \cite{Shvartsman}.
We also establish representation formulae
for these projections as in \cite{DeVoreSharpley}, where the Euclidean case is considered.

The crucial tool for the generalization of the mentioned results to $d$-sets, $n-1<d\le n$, is Markov's inequality.
We say that a set $F\subset \R^n$ {\it preserves Markov's inequality} if
the following condition holds for all positive integers $k$: for all
polynomials $p\in \mathcal{P}_k$ and all cubes $Q$ centered in $S$ with $0<\ell(Q)\le 1$,
\begin{equation}\label{markov}
\max_{F\cap Q} |\nabla p|\le c \ell(Q)^{-1} \max_{F\cap Q} |p|,
\end{equation}
where constant $c$ depends only on $F$, $n$, and $k$. 

In particular, \eqref{markov} is true when $F=\R^n$; this fact and Proposition \ref{qremez}
imply that $d$-sets, $n-1<d\le n$,
preserve Markov's inequality.
This was also shown by other methods in \cite{JonssonWallin1984}, where
function spaces on sets preserving Markov's inequality were first systematically studied.

Let $R>0$ and $Q$ be a cube which is centered in $S$ and
satisfies $\ell(Q)\le R$. 
Let also \[\{P_\beta\}_{|\beta|\le k}\]  be an
orthonormal basis of $\mathcal{P}_{k}$, $k\geq 0$, with
respect to the inner product
\[\langle f,g\rangle= \int_{Q\cap S} fg\, d\mathcal{H}^d,\qquad f,g\in \mathcal{P}_k.\]
Observe that the zero set of $p\in\mathcal{P}_{k}\setminus \{0\}$
has Hausdorff dimension at most $n-1$. Hence, the formula does
define an inner product.

Define a projection $P_{k,Q}:L^1(Q\cap S)\to \mathcal{P}_{k}$ by 
\begin{equation}\label{projs}
P_{k,Q} f :=  \sum_{|\beta|\le k} \langle f,P_\beta\rangle P_\beta
=\sum_{|\beta|\le k} \bigg(\int_{Q\cap S} fP_\beta\,d\mathcal{H}^d\bigg) P_\beta.
\end{equation}

\begin{prop}\label{pol_2}
For every $1\le u\le \infty$ and $f\in L^u(Q\cap S)$, 
$$
\bigg(\frac{1}{\mathcal{H}^d(Q\cap S)}\int_{Q\cap S}|f-P_{k,Q}(f)|^u\;d\mathcal{H}^d\bigg)^{1/u}\le C\mathcal{E}_{k+1}(f,Q)_{L^u(S)},
$$
where the constant $C$ depends on $S,R,n$ and $k$.
\end{prop}

\begin{proof}
Let us first estimate the operator norm of $P_{k,Q}$ in $L^u(Q\cap S)$. If $f\in L^u(Q\cap S)$, then
$$\Vert P_{k,Q}(f)\Vert_{L^u(Q\cap S)}\le\sum_{|\beta|\le k}\bigg|\int_{Q\cap S}fP_\beta\;d\mathcal{H}^d\bigg|\Vert P_\beta\Vert_{L^u(Q\cap S)}.$$
By H\"older's inequality
$$\Vert P_{k,Q}(f)\Vert_{L^u(Q\cap S)}\le\bigg(\sum_{|\beta|\le k}\Vert P_\beta\Vert_{L^u(Q\cap S)}\Vert P_\beta\Vert_{L^{u'}(Q\cap S)}\bigg)\Vert f\Vert_{L^u(Q\cap S)}.$$
Using \eqref{RevHolderForS-set}, orthonormality of $P_\beta$'s, and that $1/u+1/u'=1$, we get
\begin{align*}
&\Vert P_\beta\Vert_{L^u(Q\cap S)}\Vert P_\beta\Vert_{L^{u'}(Q\cap S)}\\&\le c\big((\mathcal{H}^d(Q\cap S))^{\frac{1}{u}-\frac{1}{2}}\Vert P_\beta\Vert_{L^2(Q\cap S)}\big)\big((\mathcal{H}^d(Q\cap S))^{\frac{1}{u'}-\frac{1}{2}}\Vert P_\beta\Vert_{L^2(Q\cap S)}\big)=c.
\end{align*}
It follows that $\Vert P_{k,Q}(f)\Vert_{L^u(Q\cap S)}\le c\Vert f\Vert_{L^u(Q\cap S)}$.

Let $P\in\mathcal{P}_{k}$ be such that 
$$\bigg(\frac{1}{\mathcal{H}^d(Q\cap S)}\int_{Q\cap S}|f-P|^u\;d\mathcal{H}^d\bigg)^{1/u}=\mathcal{E}_{k+1}(f,Q)_{L^u(S)}.$$
Since
$f-P_{k,Q}(f)=(f-P)-P_{k,Q}(f-P)$, it follows that
$$\bigg(\frac{1}{\mathcal{H}^d(Q\cap S)}\int_{Q\cap S}|f-P_{k,Q}(f)|^u\;d\mathcal{H}^d\bigg)^{1/u}\le(1+c) \mathcal{E}_{k+1}(f,Q)_{L^u(S)}.$$
\end{proof}

Next we  compute a representation formula for the projections. This generalizes a similar representation for projections in case of $S=\R^n$
\cite[p. 9]{DeVoreSharpley}.

\begin{lem}\label{proj_es}
There exist polynomials $h_{Q,\nu}\in\mathcal{P}_{k}$ such that 
\[
P_{k,Q} f(x) = \sum_{|\nu|\le k} \bigg(\kint_{Q\cap S} fh_{Q,\nu}\,d\mathcal{H}^d\bigg) 
\frac{(x-x_Q)^\nu}{\ell(Q)^{|\nu|}},\quad x\in\R^n,
\]
for every $f\in L^1(Q\cap S)$, and
$$\|h_{Q,\nu}\|_{L^\infty(Q)}\le C_{S,R,n,k}\, ,\, \, |\nu|\le k.$$
\end{lem}

\begin{proof}
Recall that $\{P_\beta\}_{|\beta|\le k}$  denotes an
orthonormal basis of $\mathcal{P}_{k}$. 
For every $P_\beta\in\mathcal{P}_k$, $|\beta|\le k$, by Taylor formula
\begin{equation}\label{pdf}
\begin{split}
P_\beta(x) &= \sum_{|\nu|\le k} \frac{\partial^\nu P_\beta(x_Q)}{\nu!}(x-x_Q)^\nu 
\\&=\sum_{|\nu|\le k}c_\beta^\nu (x-x_Q)^\nu.
\end{split}
\end{equation}
Let us estimate the coefficients $c_\beta^\nu$. 
Since $x_Q\in Q\subset 2Q$, by Markov-type inequality
\eqref{markov}, where $F=\R^n$, we have
\begin{align*}
|c_\beta^\nu|&\le \|\partial^\nu P_\beta\|_{L^\infty(2Q)}\le c_{n,k}\ell (2Q)^{-|\nu|}\|P_\beta\|_{L^\infty(2Q)}.
\end{align*}
Applying Proposition \ref{qremez} with exponents $r=\infty$ and $u=2$,
we get
\begin{align*}
|c_\beta^\nu|&\le c_{S,R,n,k}\ell (Q)^{-|\nu|}\bigg(\frac{1}{\mathcal{H}^d(Q\cap S)}\int_{Q\cap S} |P_\beta|^2d\mathcal{H}^d\bigg)^{1/2}\\&=c_{S,R,n,k}\ell (Q)^{-|\nu|}
\mathcal{H}^d(Q\cap S)^{-1/2}.
\end{align*}
In the last step we used that $\langle P_\beta,P_\beta\rangle=1.$
A similar computation also shows that 
\[\|P_\beta\|_{L^\infty(Q)}\le c_{S,R,n,k}\mathcal{H}^d(Q\cap S)^{-1/2}.\]

Define 
\[
h_{Q,\nu} := \sum_{|\beta|\le k} c_\beta^\nu \ell(Q)^{|\nu|} \mathcal{H}^d(Q\cap S) P_\beta\in\mathcal{P}_k,\quad |\nu|\le k.
\]
Then the previous estimates show that
$\|h_{Q,\nu}\|_{L^\infty(Q)}\le c_{S,R,n,k}$. 

Let $P_{k,Q}$ be defined as in \eqref{projs}. Then by \eqref{pdf}
\begin{align*}
P_{k,Q} f(x) &= \sum_{|\beta|\le k} \langle f,P_\beta\rangle \sum_{|\nu|\le k}c_\beta^\nu (x-x_Q)^\nu\\&= \sum_{|\nu|\le k} \bigg\langle f,\sum_{|\beta|\le k} c_\beta^\nu P_\beta\bigg\rangle (x-x_Q)^\nu\\
&= \sum_{|\nu|\le k} \bigg(\frac{1}{\mathcal{H}^d(Q\cap S)}
\int_{Q\cap S} fh_{Q,\nu}\,d\mathcal{H}^d\bigg) 
\frac{(x-x_Q)^\nu}{\ell(Q)^{|\nu|}}.
\end{align*}
\end{proof}

\subsection{A norm equivalence}\label{nors}

In this section we prove that the norm in Besov space $B_{pq}^{\alpha}(S)$ defined by \eqref{DefinitionBesovNormIntegral} is independent
of the parameter $k>\alpha$.

\begin{prop}\label{noreq}
Assume that $1\le p,q< \infty$, $1\le u\le p$ and $f\in L^p(S)$. 
Denote
\begin{align*}
N_k(f):=\bigg(\int_0^1\bigg(
\frac{\Vert\mathcal{E}_k(f,Q(\cdot,\tau))_{L^u(S)}\Vert_{L^p(S)}}{\tau^{\alpha}}\bigg)^{q}\;\frac{d\tau}{\tau}\bigg)^{\frac{1}{q}}+
\Vert f\Vert_{L^p(S)}.
\end{align*}
If $k,k'\in \N$ are such that $\alpha <\min\{k,k'\}$, then $N_k(f)$ and $N_{k'}(f)$ are
equivalent, and 
the implicit constant does not
depend on $f$.
\end{prop}

\begin{proof}
It suffices to verify the equivalence in case $k>\alpha$ and $k'=k+1$.
First of all, by Remark \ref{diskr}
\begin{align*}
E:&=\bigg(\int_0^1 \bigg(\frac{\|\mathcal{E}_{k}(f,Q(\cdot,\tau))_{L^u(S)}\|_{L^p(S)}}{\tau^\alpha}\bigg)^q\frac{d\tau}{\tau}\bigg)^{1/q}
\\&\le c\bigg\{\bigg(\sum_{j\ge 1} \big(2^{j\alpha}\|\mathcal{E}_{k}(f,Q(\cdot,2^{-j}))_{L^u(S)}\|_{L^p(S)}\big)^q\bigg)^{1/q} + \|f\|_{L^p(S)}\bigg\}.
\end{align*}
Assume that $x\in S$ and $m\ge 0$. Let $Q^m=Q(x,2^{-m})$ and $P_{k,Q^m}$ be the projections associated
to these cubes, see Section \ref{projections_sec}. By the properties of the projections
\begin{align*}
f &= f-P_{k,Q^j}(f) 
 + 
\sum_{i=0}^{j-1} \underbrace{P_{k,Q^{j-i}}\big(f-P_{{k},Q^{j-i-1}}(f)\big)}_{=:\pi_i} 
+ \underbrace {P_{k,Q^0}(f)}_{=:\pi_j}.
\end{align*}
By Lemma \ref{proj_es}, if $i=0,1,\ldots,j-1$, we have
\begin{align*}
\pi_i(y) &= 
\sum_{|\nu|=k} \kint_{Q^{j-i}\cap S} 
\big(f-P_{{k},Q^{j-i-1}}(f)\big)h_{Q^{j-i},\nu}\,d\mathcal{H}^d\,
\frac{(y-x)^\nu}{2^{k(i-j+1)}}+\rho_i(y),
\end{align*}
where $\rho_i\in \mathcal{P}_{k-1}$.
Similarly, we have
\[
\pi_j(y) = \sum_{|\nu|=k} \kint_{Q^0\cap S} 
fh_{Q^0,\nu}\,d\mathcal{H}^d\,\frac{(y-x)^\nu}{2^k}+\rho_j(y),
\]
where $\rho_j\in \mathcal{P}_{k-1}$.
Denoting $\rho=\sum_{i=0}^j \rho_i\in\mathcal{P}_{k-1}$, we have
\begin{equation}\label{sum}
\mathcal{E}_k(f,Q^j)_{L^u(S)}\le \bigg(\kint_{Q^j\cap S}
|f-\rho|^u\,d\mathcal{H}^d\bigg)^{1/u}\le  I_j(x) + II_j(x) + III_j(x).
\end{equation}
Here we have denoted
\[
I_j(x):=\bigg(\kint_{Q^j\cap S} |f-P_{k,Q^j}(f)|^u\,d\mathcal{H}^d\bigg)^{1/u},
\]
also,
\begin{align*}
II_j(x)&:=
\bigg(\kint_{Q^j\cap S} \bigg|\sum_{i=0}^{j-1}\sum_{|\nu|=k} \kint_{Q^{j-i}\cap S} 
\big(f-P_{{k},Q^{j-i-1}}(f)\big)h_{Q^{j-i},\nu}\,d\mathcal{H}^d\,\frac{(y-x)^\nu}{2^{k(i-j+1)}}\bigg|^u\,d\mathcal{H}^d(y)\bigg)^{1/u},
\end{align*}
and
\begin{align*}
III_j(x)&:=
\bigg(\kint_{Q^j\cap S}
\bigg|\sum_{|\nu|=k} \kint_{Q^0\cap S} 
fh_{Q^0,\nu}\,d\mathcal{H}^d\,\frac{(y-x)^\nu}{2^k}\bigg|^u\,d\mathcal{H}^d(y)\bigg)^{1/u}.
\end{align*}
We estimate these terms pointwise.
First, by Proposition \ref{pol_2}
\[
I_j(x) \le c_{S,n,k}\mathcal{E}_{k+1}(f,Q^j)_{L^u(S)}.
\]

Recall that the norms $||h_{Q^{j-1},\nu}||_{L^\infty(Q^{j-1})}$ are uniformly bounded
by Proposition \ref{proj_es}. Using this fact, the inequality \eqref{sset_rem} and Proposition \ref{pol_2} we get
\begin{align*}
II_j(x)&\le c\sum_{i=0}^{j-1} 2^{-ki}\kint_{Q^{j-i}\cap S} 
|f-P_{{k},Q^{j-i-1}}(f)|\,d\mathcal{H}^d\\
&\le c\sum_{i=0}^{j-1} 2^{-ki}\kint_{Q^{j-i-1}\cap S} 
|f-P_{{k},Q^{j-i-1}}(f)|\,d\mathcal{H}^d\\
&\le c\sum_{i=0}^{j-1} 2^{-ki}\mathcal{E}_{k+1}(f,Q^{j-i-1})_{L^1(S)}\le c\sum_{i=0}^{j-1} 2^{-ki}\mathcal{E}_{k+1}(f,Q^{j-i-1})_{L^u(S)}.
\end{align*}
Similar reasoning gives the estimate
\begin{align*}
III_j(x)&\le  c2^{-kj} \bigg(\frac{1}{\mathcal{H}^d(Q^0\cap S)}\int_{Q^0\cap S} |f|\,d\mathcal{H}^d\bigg).
\end{align*}
Using H\"older's inequality (if $p>1$), Fubini's theorem and \eqref{sset_rem} we get
the following estimate
\begin{equation}\label{f_est}
\begin{split}
\| III_j\|_{L^p(S)}&\le 
c2^{-kj}\bigg(\int_S \kint_{S\cap Q(x,1)} |f(y)|^p\,d\mathcal{H}^d(y)\,d\mathcal{H}^d(x)\bigg)^{1/p}\\&\le c2^{-kj}\,\| f \|_{L^p(S)}.
\end{split}
\end{equation}

By \eqref{sum}
\begin{equation*}
E\le c\bigg(\sum_{j\ge 1}(2^{j\alpha} \| I_j+II_j+III_j \|_{L^p(S)})^q\bigg)^{1/q}.
\end{equation*}
Apply now the triangle-inequality and estimate each of the
resulting three sums separately. First, we have
\begin{align*}
&\bigg(\sum_{j\ge 1}(2^{j\alpha} \| I_j(\cdot) \|_{L^p(S)})^q\bigg)^{1/q}
\\&\le c\bigg(\sum_{j\ge 1}(2^{j\alpha} \| \mathcal{E}_{k+1}(f,Q(\cdot,2^{-j}))_{L^u(S)} \|_{L^p(S)})^q\bigg)^{1/q}
\end{align*}
Next, taking into account the estimate \eqref{f_est} and the fact that $\alpha <k$, we have
\begin{align*}
\bigg(\sum_{j\ge 1}(2^{j\alpha} \| III_j(\cdot) \|_{L^p(S)})^q\bigg)^{1/q}
&\le c\|f\|_{L^p(S)}\bigg(\sum_{j\ge 1}2^{qj(\alpha-k)}\bigg)^{1/q}
\\&\le c\|f\|_{L^p(S)}.
\end{align*}
For $II_j$ we have
\begin{align*}
&\bigg(\sum_{j\ge 1}(2^{j\alpha} \| II_j(\cdot) \|_{L^p(S)})^q\bigg)^{1/q}
\\&\le c\bigg(\sum_{j\ge 1} \bigg(2^{j\alpha}\sum_{i=0}^{j-1} 2^{-ki}\|\mathcal{E}_{k+1}(f,Q(\cdot,2^{-j+i+1}))_{L^u(S)}\|_{L^p(S)}\bigg)^{q}\bigg)^{1/q}\\
&=c2^\alpha\bigg(\sum_{j\ge 1} \bigg(\sum_{i=0}^{j-1} 2^{i(\alpha-k)} 2^{\alpha(j-i-1)}\|\mathcal{E}_{k+1}(f,Q(\cdot,2^{-j+i+1}))_{L^u(S)}\|_{L^p(S)}\bigg)^{q}\bigg)^{1/q}.
\end{align*}
By the generalized Minkowski inequality 
\begin{align*}
&\bigg(\sum_{j\ge 1}(2^{j\alpha} \| II_j(\cdot) \|_{L^p(S)})^q\bigg)^{1/q}
\\&\le c\sum_{i=0}^{\infty}2^{i(\alpha-k)}\bigg(\sum_{j\ge 0} (2^{\alpha j}\|\mathcal{E}_{k+1}(f,Q(\cdot,2^{-j}))_{L^u(S)}\|_{L^p(S)})^q\bigg)^{1/q}.
\end{align*}

Collecting the estimates, we have
\begin{align*}
E&\le c\bigg\{\bigg(\sum_{j\ge 0} (2^{\alpha j}\|\mathcal{E}_{k+1}(f,Q(\cdot,2^{-j}))_{L^u(S)}\|_{L^p(S)})^q\bigg)^{1/q}+ \|f\|_{L^p(S)}\bigg\},
\end{align*}
and referring to Remark \ref{diskr} finishes the proof.
\end{proof}

\section{Restriction theorems}\label{restr_sec}
We say that a function $f\in L^1_{\rm loc}(\R^n)$ is strictly defined at $x$ if the limit
\[
\bar{f}(x)=\lim_{r\to 0+}  \kint_{Q(x,r)}f(y)\,dy=\lim_{r\to 0+} \frac{1}{|Q(x,r)|}\int_{Q(x,r)} f(y)\,dy
\]
exists.
By the Lebesgue differentiation theorem, we have $f=\bar{f}$ a.e in $\R^n$.

Let $S\subset \R^n$. At those points $x\in S$, in which $\bar{f}(x)$ exists, we define the restriction (trace) of the function $f$ on $S$ by 
$$f|_S(x):=\bar{f}(x).$$

Assume that
$f$ belongs to Besov space $B_{pq}^{\alpha}(\R^n)$
with $\alpha>0$, $1\le p<\infty$, $1\le q\le \infty$.
Then the trace $f|_S$ to a $d$-set with
$d>n-\alpha p$
is defined $\mathcal{H}^d$-almost everywhere on $S$.
In fact, the exceptional set for the Lebesgue points of $\bar f$ in $\R^n$ has zero $d$-dimensional Hausdorff measure.

This statement can be derived by using standard embeddings for these spaces \cite{Triebel1} and continuity properties of functions from Bessel potential spaces $L^\alpha_p(\R^n)$
\cite{AdamsHedberg}. 
Choose $\alpha_1\in (0,\alpha)$ such that $d>n-\alpha_1 p$
and $n/p-(\alpha-\alpha_1)>0$, then
$$
B_{pq}^{\alpha}(\R^n)\subset L_{p_1}^{\alpha_1}(\R^n),\,\,\,\,\,\,p_1=\frac{n}{n/p-(\alpha-\alpha_1)}>1.
$$
If $\alpha_1 p_1>n$, by Sobolev embedding theorem $f\in L_{p_1}^{\alpha_1}(\R^n)$ has a continuous representative in $\R^n$. For $p_1\in (1,n/\alpha_1]$ 
strictly defined functions in $L_{p_1}^{\alpha_1}(\R^n)$ have Lebesgue points $(\alpha_1,p_1)$-quasieverywhere in $\R^n$ and,
in particular, almost everywhere with respect to the $d$-dimensional Hausdorff measure \cite{AdamsHedberg}.

A similar reasoning shows that 
$\mathcal{H}^d$-almost every point, $d>n-\alpha p$, is a Lebesgue point
for a function  $\bar f$, where $f\in F_{pq}^{\alpha}(\R^n)$, $\alpha>0$,
$1<p<\infty$ and $1\le q\le \infty$.
 For more delicate continuity properties of functions from Besov and Triebel--Lizorkin spaces see \cite{Netrusov} or \cite{HedbergNetrusov}.

\subsection{Besov spaces}
In this section we prove the following restriction theorem for Besov spaces.

\begin{thm}\label{restrB}
Let $S$ be a $d$-set in $\R^n$ with $n-1<d<n$.
Let $1\le p,q<\infty$, $\alpha>(n-d)/p$ and $f\in B^{\alpha}_{pq}(\R^n)$.
Suppose also that $1\le u\le p$ and $k$ is an integer, $k>\alpha$. Then the following norm-estimate holds
\begin{equation}\label{RestrictionForBesovSpace}
\begin{split}
&\bigg(\sum_{i\ge 0} \big(2^{i(\alpha-(n-d)/p)}\|\mathcal{E}_k(f,Q(\cdot,2^{-i}))_{L^u(S)}\|_{L^p(S)}\big)^q\bigg)^{1/q}+\|f|_S\|_{L^p(S)}\\
&\le C\bigg\{\bigg(\sum_{j\ge 0} \big(2^{j \alpha}\Vert\mathcal{E}_k(f,Q(\cdot,2^{-j}))_{L^u(\R^n)}\Vert_{L^p(\R^n)}\big)^q\bigg)^{1/q}+\|f\|_{L^p(\R^n)}\bigg\},
\end{split}
\end{equation}
where $C>0$ depends on $S,n,u,p,q$ and $k$.
\end{thm}

\begin{rmk}
Theorem \ref{restrB} is essentially a special case of the result on traces of Besov spaces on $d$-sets, $0<d\le n$, proved by A. Jonsson and H. Wallin \cite{JonssonWallin1984}. However, we show that when $d>n-1$ the proof can be significantly simplified.
\end{rmk}

\noindent
The proof of Theorem \ref{restrB} uses several lemmas, in which we always assume that $S$ is a $d$-set with $n-1<d<n$.

\begin{lem}\label{EqLocalApInLpSViaLocalApInLp}
Let $f\in L^p(\R^n)$, $1\le p<\infty$, be a function such that
$\mathcal{H}^d$-every point in the $d$-set $S$ is a Lebesgue point of $f$.
Let also
 $1\le u\le p$, $d>n-\alpha p$, and $k\in\N$. Then, for every $i\in \N$, we have
\begin{equation}\label{EqLocalApInLpSViaLocalApInLp_est}
\begin{split}
&\|f|_S\|_{L^p(S)}\le c\|f\|_{L^p(\R^n)}+c\sum_{j=2}^\infty
\|\mathcal{E}_k(f,Q(\cdot,2^{-j}))_{L^u(\R^n)}\|_{L^p(S)}\,\,\,\text{and}\\
&\Vert\mathcal{E}_k(f,Q(\cdot,2^{-i-1}))_{L^u(S)}\Vert_{L^p(S)}\le
c\sum_{j=i}^\infty\Vert\mathcal{E}_k(f,Q(\cdot,2^{-j}))_{L^u(\R^n)}\Vert_{L^p(S)}.
\end{split}
\end{equation}
\end{lem}

\begin{proof}
Let   $Q_x^j$ denote a cube $Q(x,2^{-j})$ with $x\in S$ and $j\in\Z$. Let $P_{Q_x^j}:=P_{k-1,Q_x^j}$ be a
projection from $L^1(Q_x^j)$ to $\mathcal{P}_{k-1}$ such that $\mathcal{E}_k(f,Q_x^j)_{L^u(\R^n)}$
is equivalent to
$$\bigg(\kint_{Q_x^j}|f-P_{Q_x^j}f|^u\,dy\bigg)^{1/u}.$$
For the construction of these projections, see  
Proposition \ref{pol_2} or \cite{DeVoreSharpley}.
We use the following properties of these projections:
\begin{itemize}
\item[i)]$P_Q(\lambda)=\lambda\,\,\text{for any} \,\,\lambda\in \mathbb{R};$ 
\item[ii)]$|P_Q f(y)|\leq c|Q|^{-1}\int_Q |f|$ if $y\in Q$;
\item[iii)]If $Q'\subset Q$ and $|Q'|\geq c|Q|$, then for every $z\in{Q'}$
\[|P_Qf(z)-P_{Q'}f(z)|\le C\kint_{Q}|f-P_Q f|\,dy;\]
\item[iv)]$\lim_{r\to +0}P_{Q(x,r)}f(x)=f(x)$ at every Lebesgue point $x$  of $f$.
\end{itemize}
Properties i)--iii) follow from the construction of the projections. The statement iv) can then be derived from i) and ii)
as follows
\begin{equation*}
|P_{Q(x,r)}f(x)-f(x)|=|P_{Q(x,r)}[f-f(x)](x)|\leq
c\kint_{Q(x,r)}|f(y)-f(x)|dy.
\end{equation*}
The limit of the right hand side is zero at Lebesgue points of $f$.

Let $x\in S$ and $z\in Q_x^{i+1}\cap S$ be a Lebesgue point of $f$. Then by iv) 
$$
|P_{Q_x^i}f(z)-f(z)|\le |P_{Q_x^i}f(z)-P_{Q_z^{i+1}}f(z)|+\sum_{j=i+1}^\infty|P_{Q_z^{j}}f(z)-P_{Q_z^{j+1}}f(z)|.
$$
Since $z\in Q_z^{i+1}\subset Q_x^i$, the statement iii) and H\"{o}lder's inequality give
$$
|P_{Q_x^i}f(z)-P_{Q_z^{i+1}}f(z)|\le c\bigg(\kint_{Q_x^{i}}|f-P_{Q_x^i}f|^u\,dy\bigg)^{1/u}\le c\mathcal{E}_k(f,Q_x^i)_{L^u(\R^n)}.
$$
Similarly
$|P_{Q_z^{j}}f(z)-P_{Q_z^{j+1}}f(z)|\le c \mathcal{E}_k(f,Q_z^j)_{L^u(\R^n)}$ for $j\in \{i+1,\ldots\}$.
We have shown that $$
|P_{Q_x^i}f(z)-f(z)|\le c\bigg(\mathcal{E}_k(f,Q_x^i)_{L^u(\R^n)}+\sum_{j=i+1}^\infty\mathcal{E}_k(f,Q_z^j)_{L^u(\R^n)}\bigg)
$$
if $z\in Q_x^{i+1}\cap S$ is a Lebesgue point of $f$.

Since by assumption $\mathcal{H}^d$-every point in $S$ is a Lebesgue point of $f$, we can average the last inequality over the set $Q_x^{i+1}\cap S$. Then we have
\begin{align*}
&\mathcal{E}_k(f,Q_x^{i+1})_{L^u(S)}\le \bigg(\kint_{Q_x^{i+1}\cap
S}|P_{Q_x^i}f(z)-f(z)|^u\,d\mathcal{H}^d(z)\bigg)^{1/u}\\
&\le c\mathcal{E}_k(f,Q_x^i)_{L^u(\R^n)}+c\sum_{j=i+1}^\infty\bigg(\kint_{Q_x^{i+1}\cap S}\mathcal{E}_k(f,Q_z^j)^u_{L^u(\R^n)}\,d\mathcal{H}^d(z)\bigg)^{1/u}.
\end{align*}
Integrating this estimate with respect to $x$ and using H\"older's inequality gives 
\begin{align*}
&\|\mathcal{E}_k(f,Q(\cdot,2^{-i-1}))_{L^u(S)}\|_{L^p(S)}
\le
c\bigg(\int_S\mathcal{E}_k(f,Q_x^{i})_{L^u(\R^n)}^p\,d\mathcal{H}^d(x)\bigg)^{1/p}\\&\qquad\qquad+c\sum_{j=i+1}^\infty\bigg(\int_S\kint_{Q_x^{i+1}\cap
S}\mathcal{E}_k(f,Q_z^j)^p_{L^u(\R^n)}\,d\mathcal{H}^d(z)d\mathcal{H}^d(x)\bigg)^{1/p}.
\end{align*}
Using Fubini's theorem and \eqref{sset_rem}, we get
\begin{align*}
&\int_S\kint_{Q_x^{i+1}\cap
S}\mathcal{E}_k(f,Q_z^j)^p_{L^u(\R^n)}\,d\mathcal{H}^d(z)d\mathcal{H}^d(x)\\&\le c 2^{id}\int_S \int_{S} \mathcal{E}_k(f,Q_z^j)^p_{L^u(\R^n)}\chi_{Q_x^{i+1}}(z)\,d\mathcal{H}^d(x)d\mathcal{H}^d(z)\\
&\le c 2^{id}\int_{S}\mathcal{H}^d(Q_z^{i+1}\cap S)\mathcal{E}_k(f,Q_z^j)^p_{L^u(\R^n)}d\mathcal{H}^d(z)\\&\le c\int_{S}\mathcal{E}_k(f,Q_z^j)^p_{L^u(\R^n)}\,d\mathcal{H}^d(z).
\end{align*}
Collecting the estimates above, we get the second estimate in \eqref{EqLocalApInLpSViaLocalApInLp_est}.

We proceed to the verification of the remaining estimate in  \eqref{EqLocalApInLpSViaLocalApInLp_est}.
If $z\in S$ is a Lebesgue point of $f$ then, by estimating as above and using ii), we get
\begin{align*}
|f(z)|&\le |P_{Q_z^2}f(z)|+\sum_{j=2}^\infty|P_{Q_z^{j}}f(z)-P_{Q_z^{j+1}}f(z)|
\\&\le c\kint_{Q_z^2}|f(y)|\,dy+c\sum_{j=2}^\infty\mathcal{E}_k(f,Q_z^j)_{L^u(\R^n)}.
\end{align*}
By assumption, this estimate is valid $\mathcal{H}^d$-everywhere in $S$.
By integration,
\[
\|f\|_{L^p(S)}\le c\underbrace{\bigg(\int_S\kint_{Q_z^2}|f(y)|^p\,dy\,d\mathcal{H}^d(z)\bigg)^{\frac{1}{p}}}_{=:I}+c\sum_{j=2}^\infty
\|\mathcal{E}_k(f,Q(\cdot,2^{-j}))_{L^u(\R^n)}\|_{L^p(S)}.
\]
In order to finish the proof, we use Fubini's theorem for the estimate
\begin{align*}
I^p&\le c\int_S\int_{\R^n}|f(y)|^p\chi_{Q(z,1/4)}(y)\,dy\,d\mathcal{H}^d(z)\\
&\le c\int_{\R^n}|f(y)|^p\int_S\chi_{Q(z,1/4)}(y)\,d\mathcal{H}^d(z)\,dy\le c\|f\|_{L^p(\R^n)}^p.
\end{align*}
\end{proof}

\begin{lem}\label{subest}
Let $1\le u\le p<\infty$, $f\in L^1_{\mathrm{loc}}(\R^n)$, and $k\in\N$. Then for every $i\in \N$ we have
\[
\|\mathcal{E}_k(f,Q(\cdot,2^{-i}))_{L^u(\R^n)}\|_{L^p(S)}
\le c2^{i(n-d)/p}\|\mathcal{E}_k(f,Q(\cdot,2^{-i+2}))_{L^u(\R^n)}\|_{L^p(\R^n)}.
\]
Here the constant $c$ depends on $n$, $u$, $p$ and $S$.
\end{lem}

\begin{proof}
Let $\mathcal{F}=\{Q(x,2^{-i-3})\,:\,x\in S\}$.
By the $5r$-covering theorem, see e.g. \cite[p. 23]{Mat95}, there are
disjoint cubes $Q_m=Q(x_m,2^{-i-3})\in \mathcal{F}$,
$m=1,2,\ldots$ (if there are only a finite number of cubes we change the indexing), 
such that $S\subset \bigcup_{m=1}^{\infty} 5Q_m$. 
Hence
\begin{align*}
I:&=\int_{S}\mathcal{E}_k(f,Q(x,2^{-i}))_{L^u(\R^n)}^p\,d\mathcal{H}^d(x)\\&\le \sum_{m=1}^{\infty} \int_{5Q_m\cap S} \mathcal{E}_k(f,Q(x,2^{-i}))_{L^u(\R^n)}^p\,d\mathcal{H}^d(x).
\end{align*}
Notice that, if $x\in 5Q_m=Q(x_m,(5/8)2^{-i})$, then
$$Q(x,2^{-i})\subset Q(x_m,2^{-i+1})\subset Q(x,2^{-i+2}),$$ so that we have
\[
\mathcal{E}_k(f,Q(x,2^{-i}))_{L^u(\R^n)}\le c\mathcal{E}_k(f,Q(x_m,2^{-i+1}))_{L^u(\R^n)}
\le c\mathcal{E}_k(f,Q(x,2^{-i+2}))_{L^u(\R^n)}.
\]
Using the observation above and \eqref{sset_rem} we can continue as follows:
\begin{equation}\label{EstimateLemsubest}
\begin{split}
I&\le c\sum_{m=1}^{\infty} \mathcal{H}^d(5Q_m\cap S) \mathcal{E}_k(f,Q(x_m,2^{-i+1}))_{L^u(\R^n)}^p\\&\le c2^{-id}\sum_{m=1}^{\infty} \mathcal{E}_k(f,Q(x_m,2^{-i+1}))_{L^u(\R^n)}^p.\\
\end{split}
\end{equation}
From  \eqref{EstimateLemsubest} and the fact that $|Q_m|\ge c2^{-in}$ we get
\begin{align*}
I&\le c2^{i(n-d)}\sum_{m=1}^{\infty} |Q_m|\mathcal{E}_k(f,Q(x_m,2^{-i+1}))_{L^u(\R^n)}^p\\
&\le c2^{i(n-d)}\sum_{m=1}^{\infty}\int_{Q_m}\mathcal{E}_k(f,Q(x,2^{-i+2}))_{L^u(\R^n)}^p\,dx.
\end{align*}
It remains to notice that the cubes $Q_m$ are disjoint.
\end{proof}

We are ready for the proof of Theorem \ref{restrB}.
\begin{proof}
Without loss of generality we can assume that $f=\bar f$.
By Remark \ref{diskr} it is enough to estimate the sum starting from $i=3$ in \eqref{RestrictionForBesovSpace}. 
By Lemma \ref{EqLocalApInLpSViaLocalApInLp}
and Lemma \ref{subest}, for every $i\geq 3$ we have
\begin{equation*}
\Vert\mathcal{E}_k(f,Q(\cdot,2^{-i}))_{L^u(S)}\Vert_{L^p(S)}\le
c\sum_{j=i}^\infty 2^{j(n-d)/p}\Vert\mathcal{E}_k(f,Q(\cdot,2^{-j+3}))_{L^u(\R^n)}\Vert_{L^p(\R^n)}.
\end{equation*}
This results in the estimate
\begin{align*}
&\bigg(\sum_{i=3}^\infty \big(2^{i(\alpha-(n-d)/p)}\Vert\mathcal{E}_k(f,Q(\cdot,2^{-i}))_{L^u(S)}\Vert_{L^p(S)}\big)^q\bigg)^{1/q}
\\&\le c
\bigg(\sum_{i=3}^\infty 2^{i q(\alpha-(n-d)/p)}\bigg(\sum_{j=i}^\infty2^{j(n-d)/p}\Vert\mathcal{E}_k(f,Q(\cdot,2^{-j+3}))_{L^u(\R^n)}\Vert_{L^p(\R^n)}\bigg)^q\bigg)^{1/q}.
\end{align*}
Using the Hardy-type inequality \eqref{Hardy2} we can estimate
the right hand side from above by
\begin{align*}
c \bigg(\sum_{i=3}^\infty 2^{i q\alpha}\Vert\mathcal{E}_k(f,Q(\cdot,2^{-i+3}))_{L^u(\R^n)}\Vert_{L^p(\R^n)}^q\bigg)^{1/q}.
\end{align*}

It remains to estimate $\|f|_S\|_{L^p(S)}$. To this end we use Lemma \ref{EqLocalApInLpSViaLocalApInLp} and Lemma \ref{subest}. Thus, the following inequalities finish the proof
\begin{align*}
&\sum_{j=2}^\infty \|\mathcal{E}_k(f,Q(\cdot,2^{-j}))_{L^u(\R^n)}\|_{L^p(S)}\\&\le
c\sum_{j=2}^\infty2^{j(n-d)/p}\|\mathcal{E}_k(f,Q(\cdot,2^{-j+2}))_{L^u(\R^n)}\|_{L^p(\R^n)}\\
&\le c\bigg(\sum_{j=2}^\infty \big(2^{j\alpha}\|\mathcal{E}_k(f,Q(\cdot,2^{-j+2}))_{L^u(\R^n)}\|_{L^p(\R^n)}\big)^q\bigg)^{1/q}\bigg(\sum_{j=2}^\infty 2^{jq'((n-d)/p-\alpha)}\bigg)^{1/q'}.
\end{align*}
Note, that since $\alpha>(n-d)/p$, the last series converges.
\end{proof}

\subsection{Triebel--Lizorkin spaces}
Here we prove the restriction theorem for Triebel--Lizorkin spaces.
Our method of proof is different from that of \cite{SakaI} and \cite{Mamedov}, where the jet technique from \cite{JonssonWallin1984} is used.

\begin{thm}\label{ctrace}
Let $S\subset\R^n$ be a $d$-set with 
$n-1<d<n$.
Let $1< p<\infty$, $1\le q\le \infty$,  $\alpha>(n-d)/p$ and $k>\alpha$ be an integer. Then, for every
$f\in F^{\alpha}_{pq}(\R^n)$, we have
\begin{equation}\label{numbest}
\bigg(\sum_{i=0}^\infty \big(2^{i(\alpha-(n-d)/p)}\Vert\mathcal{E}_k(f,Q(\cdot,2^{-i}))_{L^1(S)}\Vert_{L^p(S)}\big)^p\bigg)^{1/p} + \|f|_S\|_{L^p(S)}
\le c\|f\|_{F^{\alpha}_{pq}(\R^n)}.
\end{equation} 
Here the constant
$c$ depends on $n$, $S$, $p$ and $\alpha$.
\end{thm}

Let $S$ be a $d$-set with $n-1<d<n$. Since $S$ is a closed set, its complement admits a Whitney decomposition, see e.g. \cite{Stein}.
That is, there is a family  $\mathcal{W}_S$ of  dyadic cubes
whose interiors are pairwise disjoint and 
$\mathbb{R}^n\setminus S=\bigcup_{Q\in \mathcal{W}_S}Q$. Furthermore,
if $Q\in \mathcal{W}_S$, then
\begin{equation}\label{dist_est}
{\rm diam}\,Q\le {\rm dist}(Q,S)\le 4\,{\rm diam}\,Q.
\end{equation}

We write $\ell(Q)\overset{\kappa}{\sim} 2^{-i}$ if
$Q$ is a cube in $\R^n$ and $2^{-i}/5\kappa\le \mathrm{diam}\,Q\le 2^{-i}$.
The lemma below uses the fact that $S$ is $\kappa$-porous, see
Remark \ref{porous_lem}.

\begin{lem}\label{cop}
Assume that $x\in S$ and $i\in\N$. Then there is $Q\in\mathcal{W}_S$ such that
$\ell(Q)\overset{\kappa}{\sim} 2^{-(i+1)}$ and 
$Q\subset Q(x,2^{-i})$. 
\end{lem}

\begin{proof}
By $\kappa$-porosity of $S$, 
there is $y\in Q(x,2^{-i-1})$ such that
$Q(y,2^{-i-1}/\kappa)\subset \R^n\setminus S$.
Let $Q\in\mathcal{W}_S$ be a cube containing the point $y$. Then
\[
2^{-i-1}/{\kappa}-\mathrm{diam}\,Q\le \mathrm{dist}(y,S)-\mathrm{diam}\,Q\le \mathrm{dist}(Q,S).
\]
On the other hand,
\[
\mathrm{dist}(Q,S)\le \mathrm{dist}(y,S)\le \|x-y\|_\infty\le 2^{-i-1}.
\]
By \eqref{dist_est}, 
$2^{-i-1}/5\kappa\le \mathrm{diam}\,Q\le 2^{-i-1}$, that is,
$\ell(Q)\overset{\kappa}{\sim} 2^{-(i+1)}$. 

It is also easy to see that $Q\subset Q(x,2^{-i})$.
\end{proof}

For $f\in L^1_{\rm loc}(\R^n)$ and $\alpha>0$ denote $$f^\sharp_\alpha(x)=\sup_{r>0} r^{-\alpha}\mathcal{E}_k(f,Q(x,r))_{L^1(\R^n)},\quad x\in\R^n,\,k=[\alpha]+1.$$

\begin{lem}\label{subestTL}
Let $1< p<\infty$, $\alpha>0$, $k=[\alpha]+1$ and $f\in L^1_{\rm loc}(\R^n)$. Then for every $i\in \N$,
\begin{align*}
\|\mathcal{E}_k(f,Q(\cdot,2^{-i}))_{L^1(\R^n)}\|_{L^p(S)}
\le c2^{-i(\alpha -(n-d)/p)}\bigg(\int_{\cup \{Q\in\mathcal{W}_S\,:\, \ell(Q)\overset{\kappa}\sim 2^{-i-4}\}}f^\sharp_\alpha(x)^pdx\bigg)^{1/p},
\end{align*}
where the constant $c$ depends on $\alpha$, $p$, $n$ and $S$.
\end{lem}

\begin{proof}
Repeating the proof of Lemma \ref{subest} we arrive at (\ref{EstimateLemsubest}), i.e.
\[
I:=\int_{S}\mathcal{E}_k(f,Q(x,2^{-i}))_{L^1(\R^n)}^pd\mathcal{H}^d(x)\le c2^{-id}\sum_{m=1}^{\infty} \mathcal{E}_k(f,Q(x_m,2^{-i+1}))_{L^1(\R^n)}^p,
\]
where $Q_m=Q(x_m, 2^{-i-3})$, $m=1,2,\dots$, are disjoint cubes centered at $S$.

By Lemma \ref{cop}, for every $Q_m$, there is  $R_m\in\mathcal{W}_S$ such that
$\ell(R_m)\overset{\kappa}{\sim} 2^{-i-4}$ and 
$R_m\subset Q_m$. 
If $x\in R_m$, then $Q(x_m,2^{-i+1})\subset Q(x,2^{-i+2})$ and  therefore
\begin{equation}\label{p3}
\mathcal{E}_k(f,Q(x_m,2^{-i+1}))_{L^1(\R^n)}\le c \mathcal{E}_k(f,Q(x,2^{-i+2}))_{L^1(\R^n)}\le c2^{-i\alpha}f^\sharp_\alpha(x).
\end{equation}
Since $\ell(R_m)\overset{\kappa}{\sim} 2^{-i-4}$, we have $|R_m|\ge c2^{-in}$. By using this and (\ref{p3}), we get
\begin{align*}
I&\le c2^{i(n-d)}\sum_{m=1}^{\infty} |R_m|\mathcal{E}_k(f,Q(x_m,2^{-i+1}))_{L^1(\R^n)}^p\\
&\le c2^{i(n-d-\alpha p)}\sum_{m=1}^{\infty}\int_{R_m}f^\sharp_\alpha(x)^p\,dx\\
&\le c2^{i(n-d-\alpha p)}\int_{\cup \{Q\in\mathcal{W}_S\,:\, \ell(Q)\overset{\kappa}\sim 2^{-(i+4)}\}}f^\sharp_\alpha(x)^pdx.
\end{align*}
Taking the $p$'th roots yield the required estimate.
\end{proof}

We are now ready for the proof of Theorem \ref{ctrace}.

\begin{proof}
Let $f\in F^{\alpha}_{pq}(\R^n)$.
We assume that $f=\bar f$.
Since
\[
\mathcal{E}_k(f,Q)_{L^1(S)}\le \mathcal{E}_{[\alpha]+1}(f,Q)_{L^1(S)},\qquad k>[\alpha]+1,
\]
we can also assume that $k=[\alpha]+1$.

Using the boundedness of the Hardy--Littlewood maximal operator for $p>1$ and the definition of Triebel--Lizorkin spaces, see Section \ref{SecFunctionSpacesRn}, we obtain
\[
\|f^\sharp_\alpha\|_{L^p(\R^n)}+\|f\|_{L^p(\R^n)}\le 
c\{\|g\|_{L^p(\R^n)} +\|f\|_{L^p(\R^n)}\big\}\le 
c \|f\|_{F^{\alpha}_{p\infty}}\le c\|f\|_{F^{\alpha}_{pq}}.
\]
Here $g(x)=\sup\{t^{-\alpha}\mathcal{E}_{[\alpha]+1}(f,Q(x,t))_{L^1(\R^n)}:\, 0<t\le 1\}$.

Hence, it suffices to verify that the left hand side of
\eqref{numbest} is bounded (up to a constant) by 
$$\|f_\alpha^\sharp\|_{L^p(\R^n)} + \|f\|_{L^p(\R^n)}.$$
Also, by Remark \ref{diskr} it is enough to estimate the sum starting from $i=2$ in \eqref{numbest}.

By Lemma \ref{EqLocalApInLpSViaLocalApInLp} and Lemma \ref{subestTL}, for every
$i\in\N$,
\begin{align*}
&\Vert\mathcal{E}_k(f,Q(\cdot,2^{-i-1}))_{L^1(S)}\Vert_{L^p(S)}\le
c\sum_{j=i}^\infty\Vert\mathcal{E}_k(f,Q(\cdot,2^{-j}))_{L^1(\R^n)}\Vert_{L^p(S)}\\
&\le c\sum_{j=i}^\infty2^{-j(\alpha -(n-d)/p)}\bigg(\int_{\cup \{Q\in\mathcal{W}_S\,:\, \ell(Q)\overset{\kappa}\sim 2^{-j-4}\}}f^\sharp_\alpha(x)^pdx\bigg)^{1/p}.
\end{align*}
By Hardy-type inequality (see Lemma \ref{LemmaHardyInequality}) with $\sigma=\alpha p-(n-d)>0$,
\begin{equation}\label{dsn}
\begin{split}
&\sum_{i=1}^\infty \big(2^{i(\alpha-(n-d)/p)}\Vert\mathcal{E}_k(f,Q(\cdot,2^{-i-1}))_{L^1(S)}\Vert_{L^p(S)}\big)^p
\\&\le c\sum_{i=1}^\infty 2^{i(\alpha p-(n-d))}\bigg(\sum_{j=i}^\infty2^{-j(\alpha -(n-d)/p)}\bigg(\int_{\cup \{Q\in\mathcal{W}_S\,:\, \ell(Q)\overset{\kappa}\sim 2^{-j-4}\}}f^\sharp_\alpha(x)^pdx\bigg)^{1/p}\bigg)^p
\\&\le c \sum_{i=1}^\infty\int_{\cup \{Q\in\mathcal{W}_S\,:\, \ell(Q)\overset{\kappa}\sim 2^{-i-4}\}}f^\sharp_\alpha(x)^pdx.
\end{split}
\end{equation}

Let us
denote $U_i:=\cup \{\mathrm{int}\,Q\,:\,Q\in\mathcal{W}_S\text{ and }\ell(Q)\overset{\kappa}\sim 2^{-i-4}\}$, and
choose $k_0\in \N$ such that $2^{-k_0}<1/5\kappa$.
Then we claim that
\begin{equation}\label{dis}
U_i\cap U_{i'}=\emptyset,\qquad \text{ if }i\not=i'\text{ and }i,i'\equiv k\,\mathrm{mod}\,k_0.
\end{equation}
To verify this claim, let $i>i'$ be such that $i,i'\equiv k\,\mathrm{mod}\,k_0$. Then $i-i'\ge k_0$.
In particular, if $Q\in \mathcal{W}_S$, $\ell(Q)\overset{\kappa}{\sim} 2^{-i-4}$, and 
$Q'\in \mathcal{W}_S$, $\ell(Q')\overset{\kappa}{\sim} 2^{-i'-4}$, then
\[
\mathrm{diam}\,Q\le 2^{-i-4}\le 2^{-k_0}2^{-i'-4}<2^{-i'-4}/5\kappa\le \mathrm{diam}(Q').
\]
It follows that $Q\not=Q'$. Since the interiors of Whitney cubes are pairwise disjoint, we find that 
$\mathrm{int}\,Q\cap \mathrm{int}\,Q'=\emptyset$. Hence, \eqref{dis} holds.

From \eqref{dis} it follows that
\begin{equation}\label{dsm}
\begin{split}
\sum_{i=1}^\infty \int_{\cup \{Q\in\mathcal{W}_S\,:\, \ell(Q)\overset{\kappa}\sim 2^{-i-4}\}}f^\sharp_\alpha(x)^pdx
&=\sum_{k=0}^{k_0-1}  \sum_{\substack{i\ge 1\\ i\equiv k\,\mathrm{mod}\,k_0}}\int_{U_i}f^\sharp_\alpha(x)^pdx
\\&\le k_0\|f_\alpha^\sharp\|_{p}^p.
\end{split}
\end{equation}

Combining the estimates \eqref{dsn} and \eqref{dsm}, we find that
\[
\bigg(\sum_{i=2}^\infty \big(2^{i(\alpha-(n-d)/p)}\|\mathcal{E}_k(f,Q(\cdot,2^{-i}))_{L^1(S)}\|_{L^p(S)}\big)^p\bigg)^{1/p}\le c\|f_\alpha^\sharp\|_p.
\]
Next we estimate $\|f|_S\|_{L^p(S)}$. By Lemma \ref{EqLocalApInLpSViaLocalApInLp},
\[
\|f|_S\|_{L^p(S)}\le c\|f\|_{L^p(\R^n)}+c\underbrace{\sum_{j=2}^\infty
\|\mathcal{E}_k(f,Q(\cdot,2^{-j}))_{L^1(\R^n)}\|_{L^p(S)}}_{=:II}.
\]
The second term is estimated using Lemma \ref{subestTL} and \eqref{dsm} as follows:
\begin{align*}
II&\le c\sum_{j=2}^\infty 2^{-j(\alpha-(n-d)/p)}\bigg(\int_{\cup \{Q\in\mathcal{W}_S\,:\, \ell(Q)\overset{\kappa}\sim 2^{-j-4}\}}f^\sharp_\alpha(x)^pdx\bigg)^{1/p}\\
&\le c\bigg(\sum_{j=2}^\infty\int_{\cup \{Q\in\mathcal{W}_S\,:\, \ell(Q)\overset{\kappa}\sim 2^{-j-4}\}}f^\sharp_\alpha(x)^pdx\bigg)^{1/p}\bigg(\sum_{j=2}^\infty 2^{-jp'(\alpha-(n-d)/p)}\bigg)^{1/p'} \\
&\le c\|f^\sharp_\alpha\|_{L^p(\R^n)}.
\end{align*}
In the last step we also used the estimate $\alpha -(n-d)/p>0$.
\end{proof}

\section{Extension operator}\label{ext_sec}

\subsection{Constructing the extension operator}

The construction of the extension operator is based on a modification of the 
Whitney extension method. A similar modification was first suggested by Shvartsman in \cite{Shvartsman78}
and developed later in \cite{Shvartsman}.
Our work extends the results in the latter paper to the case
of $d$-sets, $n-1<d<n$.

Recall that $\mathcal{W}_S$ denotes a Whitney decomposition 
of $\mathbb{R}^n\setminus S$, so that $$\R^n\setminus S = \bigcup_{Q\in\mathcal{W}_S} Q,$$
and if $Q\in \mathcal{W}_S$, then
\[{\rm diam}\,Q\le {\rm dist}(Q,S)\le 4\,{\rm diam}\,Q.\]
Moreover, there is a constant $N>0$ such that 
every point of $\mathbb{R}^n\setminus S$ is covered by 
at most $N$ cubes in $\mathcal{W}_S$.

Let $\Phi:=\{\varphi_Q:Q\in \mathcal{W}_S\}$ be a smooth partition of unity
which is subordinate to the Whitney decomposition $\mathcal{W}_S$.
Then, in particular, 
\[\chi_{\R^n\setminus S}= \sum_{Q\in\mathcal{W}_S} \varphi_Q
\]
and
$\mathrm{supp}\,\varphi_Q\subset Q^*=\frac{9}{8}Q$ for every
$Q\in\mathcal{W}_S$.

To every $Q=Q(x_Q,r_Q)\in \mathcal{W}_S$, assign the cube
$
a(Q):=Q(a_Q,r_Q/2),
$
where $a_Q\in S$ is such that
$\|x_Q-a_Q\|_\infty =\mathrm{dist}(x_Q,S)$.

Let $\Delta>0$ be a parameter which will be chosen later. Then 
\begin{equation}\label{ii)}
\mathcal{H}^d(a(Q)\cap S)\geq 
c\, r_Q^d,\qquad \text{if }{\rm diam}\,Q\le\Delta,
\end{equation}
and 
\begin{equation}\label{i)}
a(Q)\subset \mathrm{int}(10Q),\qquad  Q\in\mathcal{W}_S.
\end{equation}
Indeed, \eqref{ii)}   follows from \eqref{sset_rem}, 
and the constant $c>0$ depends on  $\Delta,S$ and $n$.
Furthermore, if 
$y\in a(Q)$, then \begin{align*}
\|y-x_Q\|_\infty&\le \|y-a_Q\|_\infty+\|a_Q-x_Q\|_\infty\\&< r_Q+{\rm dist}(x_Q,S)\le 2r_Q+{\rm dist}(Q,S),
\end{align*}
and, since
$\mathrm{dist}(Q,S)\le 4 \mathrm{diam} Q=8 r_Q$,  \eqref{i)} is true.

For every $Q\in\mathcal{W}_S$ and $k\in\N_0$, define an operator $P_{k,Q}:L^1(a(Q)\cap S)\to \mathcal{P}_{k}$ by
\[
P_{k,Q}:=
\begin{cases}
P_{k,a(Q)},\quad &\text{if }{\rm diam}\,Q\le \Delta;\\
0,\quad &\text{if }{\rm diam}\,Q>\Delta.
\end{cases}
\] 
Here the projection $P_{k,a(Q)}:L^1(a(Q)\cap S)\to \mathcal{P}_{k}$,
which is
associated with a cube $a(Q)$ centered in $S$, is
constructed in Section \ref{projections_sec}.

If $f\in L^1_{\mathrm{loc}}(S)$
and $k\in\N$, we define \begin{equation}\label{DefExtensionOperator}
{\rm Ext}_{k,S}f(x):=\begin{cases}
f(x), \quad &\text{if }x\in S;\\
\sum_{Q\in \mathcal{W}_S}\varphi_Q(x)P_{k-1,Q}f(x), \quad &\text{if }x\in\mathbb{R}^n\setminus S.
\end{cases}
\end{equation}
Observe that  \eqref{DefExtensionOperator}
induces a linear extension operator ${\rm Ext}_{k,S}$.
Also, from the properties of dyadic cubes it follows that
\begin{equation}\label{supp}
\mathrm{supp}({\rm Ext}_{k,S}f)\subset \bigcup_{x\in S}Q(x,8\Delta).
\end{equation}

\subsection{Local approximations of the extension}

In this section we estimate the local approximations of the extension 
\[
\tilde{f}:={\rm Ext}_{k,S}f,
\]
where $k\in\N$ and  $f\in L^1_{\mathrm{loc}}(S)$.
To this end, we follow the proof of Theorem 3.6 in \cite{Shvartsman} quite closely. But since some of the required modifications are non-trivial,
we present detailed proofs nevertheless.

Let $x\in\mathbb{R}^n$ and $t>0$, denote by $a_x\in S$ a point
for which ${\rm dist}(x,S)=\|x-a_x\|_\infty$.
Define
$$
K^{(x,t)}:=Q(a_x,r^{(x,t)}),\qquad r^{(x,t)}:=50\max(80t,{\rm dist}(x,S)).
$$

One of the main results in this section is the following.

\begin{prop}\label{th36}
Let $f\in L^u_{\mathrm{loc}}(S)$, $1\le u< \infty$, and  $k\in\N$. 
Then
\begin{equation}\label{rest}
\tilde f|_S = f\qquad \mathcal{H}^d\text{-a.e. in }S.
\end{equation}
Furthermore, 
if $x\in\mathbb{R}^n$ and $t>0$, then
\begin{equation}\label{norm}
\begin{split}
&\mathcal{E}_k(\tilde{f},Q(x,t))_{L^u(\R^n)}\\&\le c \frac{t^k}{t^k+{\rm dist}(x,S)^k}\left\{\begin{array}{l}
\mathcal{E}_k(f,K^{(x,t)})_{L^u(S)},\,\,r^{(x,t)}\in (0,\Delta]; \\
\\
\mathcal{E}_0(f,K^{(x,t)})_{L^u(S)},\,\,r^{(x,t)}\in (\Delta,1024\Delta].
\end{array}
\right.
\end{split}
\end{equation}
Here the constant $c$ depends on $S,n,u,\Delta$ and $k$.
\end{prop}

The proof of Proposition \ref{th36} proceeds in various
cases, which are treated throughout this section. At the end of the section
we collect all of these separate cases.

To begin with, for a given cube $K$, define some families of Whitney cubes:
\begin{align*}
\mathcal{Q}_1(K)&:=\{Q\in \mathcal{W}_S\,:\,Q\cap K\neq\emptyset\},\\
\mathcal{Q}_2(K)&:=\{Q\in \mathcal{W}_S\,:\,\exists\, Q'\in \mathcal{Q}_1(K)\, \text{such that}\,Q'\cap Q^*\neq\emptyset\},\\
\mathcal{Q}_3(K)&:=\{Q\in \mathcal{Q}_2(K)\,:\,\mathrm{diam}\,Q\le\Delta\}.
\end{align*}

The proof of the following lemma
can be found in \cite[Lemma 3.7]{Shvartsman}.

\begin{lem}\label{kest}
Let $K$ be a cube centered in $S$. Then, for every $Q\in\mathcal{Q}_2(K)$,
\[\mathrm{diam}\,Q\le 2\,\mathrm{diam}\,K,\qquad \Vert x_K-x_Q\Vert_\infty\le \frac{5}{2}\mathrm{diam}\,K.\]
\end{lem}

\begin{lem}\label{pointsp}
Let $K$ be a cube centered in $S$. Then for
every $f\in L^u_{\mathrm{loc}}(S)$, $1\le u<\infty$, and every polynomial $P_0$,
\[
\sum_{Q\in\mathcal{Q}_3(K)} \frac{|Q|}{\mathcal{H}^d(a(Q)\cap S)}\|f-P_0\|_{L^u(a(Q)\cap S)}^u \le cr_K^{n-d} \|f-P_0\|_{L^u(25K\cap S)}^u ,
\]
where the constant $c$ depends on $S,\Delta$ and $n$.
\end{lem}

\begin{proof}
Suppose that $Q\in\mathcal{Q}_3(K)\subset \mathcal{Q}_2(K)$. If $y\in 10Q$
then, by Lemma \ref{kest}, we have
\begin{align*}
\|y-x_K\|_\infty&\le \mathrm{diam}(10Q)/2 + (5/2)\mathrm{diam}\,K
\\&\le 10\mathrm{diam}\,K+(5/2)\mathrm{diam}\,K=(25/2)\mathrm{diam}\,K=25r_K.
\end{align*}
Hence, $10Q\subset Q(x_K,25r_K)=25K$. Since, by \eqref{i)}, $a(Q)\cap S\subset 10Q\cap S$, we also have
\[
a(Q)\cap S\subset 25K\cap S.
\]
Recall that $\mathrm{diam}\,Q\le \Delta$ if $Q\in\mathcal{Q}_3(K)$. Thus, by \eqref{ii)},
$\mathcal{H}^d(a(Q)\cap S) \ge cr_Q^d$.

By the previous observations, we obtain
\begin{align*}
&\sum_{Q\in\mathcal{Q}_3(K)} \frac{|Q|}{\mathcal{H}^d(a(Q)\cap S)}\|f-P_0\|_{L^u(a(Q)\cap S)}^u\\
&\le c_{S,n,\Delta}\int_{25K\cap S} \sum_{Q\in\mathcal{Q}_3(K)} r_Q^{n-d} \chi_{a(Q)\cap S}(x)|f(x)-P_0(x)|^u\,d\mathcal{H}^d(x)\\
&\le c_{S,n,\Delta}\bigg\|\sum_{Q\in\mathcal{Q}_3(K)} r_Q^{n-d}\chi_{a(Q)\cap S}\bigg\|_{L^\infty(S)}
\cdot\int_{25K\cap S} |f(x)-P_0(x)|^ud\mathcal{H}^d(x).
\end{align*}
Denote
$$I(x):=\sum_{Q\in\mathcal{Q}_3(K)} r_Q^{n-d}\chi_{a(Q)\cap S}(x),\quad x\in S.$$
To finish the proof of the lemma it remains to show that $\|I\|_{L^\infty(S)}\le cr_K^{n-d}$. To this end write $I(x)$ in the following way
\begin{align*}
I(x)=
2^{d-n}\sum_{j=-\infty}^\infty  \sum_{Q\in \mathcal{D}_j\cap \mathcal{Q}_3(K)} 2^{j(d-n)} \chi_{a(Q)\cap S}(x),
\end{align*}
where $\mathcal{D}_j=\{Q\in\mathcal{D}\,:\,\ell(Q)=2^{-j}\}$.

Choose $j_0\in \Z$ such that
$2^{-j_0-1}< \mathrm{diam}\,K\le 2^{-j_0}$. If $Q\in\mathcal{D}_j\cap \mathcal{Q}_3(K)$, then, by Lemma \ref{kest},
$2^{-j}=\mathrm{diam}\,Q\le 2\mathrm{diam}\,K\le 2^{-j_0+1}$. Hence, in fact,
\begin{equation}\label{I(x)}
I(x)=2^{d-n}\sum_{j=j_0-1}^\infty 2^{j(d-n)}  \sum_{Q\in \mathcal{D}_j\cap \mathcal{Q}_3(K)} 
\chi_{a(Q)\cap S}(x).
\end{equation}

Let us also note that, for any $x\in S$ and $j\geq j_0-1$,
\begin{equation}\label{eqSumOfChi}
\sum_{Q\in \mathcal{D}_j\cap \mathcal{Q}_3(K)}
\chi_{a(Q)\cap S}(x)= {\rm card}\{Q\in \mathcal{D}_j\cap \mathcal{Q}_3(K):\,x\in a(Q)\cap S\}\le 20^n.
\end{equation}
Indeed, the upper bound here can be found in the following way. Suppose that $Q\in \mathcal{D}_j\cap \mathcal{Q}_3(K)$ is such that $x\in a(Q)\cap S$. Then $x\in 10Q$ and, if $y\in Q\subset 10Q$, we have
$\|y-x\|_\infty\le \mathrm{diam}(10Q) = 10 \cdot2^{-j}$. Hence, it follows that $Q\subset Q(x,10\cdot 2^{-j})$.
Since the side length of the cubes in $\mathcal{D}_j\cap \mathcal{Q}_3(K)$ is $2^{-j}$
and the interiors of these cubes are disjoint, there exist at most $20^n$ such cubes which are contained in a cube $Q(x,10\cdot 2^{-j}) $.

By \eqref{I(x)} and \eqref{eqSumOfChi}
$$
\|I\|_{L^\infty(S)}\le c_{d,n}2^{j_0(d-n)}<c_{d,n}(\mathrm{diam}\,K)^{n-d}=c_{d,n}r_K^{n-d}.
$$

\end{proof}

\begin{lem}\label{knine}
Assume that $f\in L^u_{\mathrm{loc}}(S)$, $1\le u <\infty$, and $k\in\N$. Then for every cube $K=Q(x_K,r_K)$ centered in $S$ and every polynomial $P_0\in  \mathcal{P}_{k-1}$
$$\sum_{Q\in\mathcal{Q}_3(K)}\Vert P_{k-1,Q}f-P_0\Vert^u_{L^u(Q)}\le cr_K^{n-d}\Vert f-P_0\Vert^u_{L^u(25K\cap S,\mathcal{H}^d)}.$$
Here the constant $c$ depends on $S,\Delta,n,u$ and $k$.
\end{lem}

\begin{proof}
Denote $p:=P_{k-1,Q}f-P_0\in\mathcal{P}_{k-1}$.
Let $Q\in \mathcal{Q}_3(K)\subset \mathcal{W}_S$. Then
$\mathrm{diam}\,Q\le \Delta$ and therefore $P_{k-1,Q}=P_{k-1,a(Q)}$.
By \eqref{i)} and 
Proposition \ref{qremez},
\begin{align*}
\frac{1}{|Q|}\int_{Q}|p|^u\;dx
\le
\frac{10^n}{|10Q|}\int_{10Q}|p|^u\;dx
\le 
\frac{c}{\mathcal{H}^d(a(Q)\cap S)}\int_{a(Q)\cap S}|p|^u\;d\mathcal{H}^d.
\end{align*}
By Proposition \ref{pol_2}
\begin{align*}
\Vert p\Vert_{L^u(a(Q)\cap S)}&\le \Vert P_{k-1,a(Q)}f-f\Vert_{L^u(a(Q)\cap S)}+\Vert f-P_0\Vert_{L^u(a(Q)\cap S)}\\&\le c\mathcal{H}^d(a(Q)\cap S)^{1/u}
\mathcal{E}_k(f,a(Q))_{L^u(S)}+\Vert f-P_0\Vert_{L^u(a(Q)\cap S)}\\
&\le c\Vert f-P_0\Vert_{L^u(a(Q)\cap S)}.
\end{align*}
Hence, for every $Q\in\mathcal{Q}_3(K)$
we have
$$\frac{1}{|Q|}\int_{Q}|p|^u\;dx\le \frac{c}{\mathcal{H}^d(a(Q)\cap S)}\int_{a(Q)\cap S}|f-P_0|^u\;d\mathcal{H}^d.$$
Summing these estimates gives
$$\sum_{Q\in\mathcal{Q}_3(K)}\Vert p\Vert^u_{L^u(Q)}\le c\sum_{Q\in\mathcal{Q}_3(K)}\frac{|Q|}{\mathcal{H}^d(a(Q)\cap S)}\Vert f-P_0\Vert^u_{L^u(a(Q)\cap S)}.$$
To finish the proof, it remains to apply Lemma \ref{pointsp}.
\end{proof}

\begin{lem}\label{chainest}
Let $f\in L^u_{\mathrm{loc}}(S)$ and $1\le u <\infty$. Then for every cube $K$ and every polynomial $P_0\in  \mathcal{P}_{k-1}$,
$$\Vert \tilde{f}-P_0\Vert^u_{L^u(K\setminus S)}\le c_{n,k,u} \sum_{Q\in\mathcal{Q}_2(K)}\Vert P_{k-1,Q}f-P_0\Vert^u_{L^u(Q)}.$$
\end{lem}

The proof of this lemma can be found  in \cite[p. 1222]{Shvartsman}.

\begin{lem}\label{Proposition310}
Let $f\in L^u_{\mathrm{loc}}(S)$, $1\le u <\infty$. Then for every cube $K$, centered in $S$, 
with $\mathrm{diam}\, K\le \Delta/2$, and every polynomial $P_0\in  \mathcal{P}_{k-1}$,
\[
\kint_{K}|\tilde{f}-P_0|^udx \le c \kint_{25K\cap S} |f-P_0|^u d\mathcal{H}^d,
\]
where the constant $c$ depends on $S,\Delta,n,u$ and $k$.
\end{lem}

\begin{proof}
From the fact that $|S|=0$ and Lemma \ref{chainest} we obtain
$$
\Vert \tilde{f}-P_0\Vert^u_{L^u(K)}=
\Vert \tilde{f}-P_0\Vert^u_{L^u(K\setminus S)}\le c \sum_{Q\in\mathcal{Q}_2(K)}\Vert P_{k-1,Q}f-P_0\Vert^u_{L^u(Q)}.$$
Let $Q\in\mathcal{Q}_2(K)$. 
By Lemma \ref{kest}, 
$\mathrm{diam}\,Q\le 2\mathrm{diam} K\le \Delta$, so, in fact, $\mathcal{Q}_2(K)=\mathcal{Q}_3(K)$. Hence, by Lemma \ref{knine},
\begin{align*}
&\sum_{Q\in\mathcal{Q}_2(K)}\Vert P_{k-1,Q}f-P_0\Vert^u_{L^u(Q)}\\&=\sum_{Q\in\mathcal{Q}_3(K)}\Vert P_{k-1,Q}f-P_0\Vert^u_{L^u(Q)}\le cr_K^{n-d}\Vert f-P_0\Vert^u_{L^u(25K\cap S)}.
\end{align*}
It remains to use the estimate \eqref{sset_rem}.
\end{proof}

\begin{lem}\label{th311}
Let $f\in L^u_{\mathrm{loc}}(S)$ and $1\le u< \infty$. Then for every cube $K=Q(x_K,r_K)$, centered in $S$, with $\mathrm{diam}\,K\le \Delta/2$
$$\mathcal{E}_k(\tilde{f},K)_{L^u(\R^n)}\le c \mathcal{E}_k(f,25K)_{L^u(S)}.$$
Here the constant $c$ depends on $S,\Delta,n,u$  and $k$.
\end{lem}

\begin{proof}
Let $P_0\in \mathcal{P}_{k-1}$ be such that
$$
\mathcal{E}_k(f,25K)_{L^u(S)}=\bigg(\frac{1}{\mathcal{H}^d(25K\cap S)}\int_{25K\cap S}|f-P_0|^u\;d\mathcal{H}^d\bigg)^{1/u}.
$$
Then Lemma \ref{Proposition310} and 
\eqref{sset_rem} imply
\begin{align*}
\mathcal{E}_k(\tilde{f},K)_{L^u(\R^n)}&\le 
c_{n,u}r_K^{-\frac{n}{u}}\Vert \tilde{f}-P_0\Vert_{L^u(K)}\\&\le c r_K^{-\frac{d}{u}}\Vert f-P_0\Vert_{L^u(25K\cap S)}\le c\mathcal{E}_k(f,25K)_{L^u(S)}.
\end{align*}
\end{proof}

\begin{lem}\label{p312}
Let $f\in L^u_{\mathrm{loc}}(S)$ and $1\le u <\infty$. Then for every cube $K$,
centered in $S$ with
$\mathrm{diam}(K)\le 1024\Delta$,
$$
\kint_K|\tilde{f}|^{u}\,dx\le 
c\kint_{25K\cap S}
|f|^{u}\,d\mathcal{H}^d,
$$
where $\tilde f=\mathrm{Ext}_{k,S} f$ and
the constant $c$ depends on $S,\Delta,n,u$ and $k$.
\end{lem}

\begin{proof}
From the fact that $|S|=0$ and
Lemma \ref{chainest} with $P_0=0$ we get
$$
\Vert \tilde{f}\Vert^u_{L^u(K)}=\Vert \tilde{f}\Vert^u_{L^u(K\setminus S)}\le c \sum_{Q\in\mathcal{Q}_2(K)}\Vert P_{k-1,Q}f\Vert^u_{L^u(Q)}.$$
Taking into account that $P_{k-1,Q}f=0$ if ${\rm diam}\,Q>\Delta$, $Q\in\mathcal{Q}_2(K)$, and the statement of Lemma \ref{knine}, we have
$$\Vert \tilde{f}\Vert^u_{L^u(K,\;dx)}=\Vert \tilde{f}\Vert^u_{L^u(K\setminus S)}\le c \sum_{Q\in\mathcal{Q}_3(K)}\Vert P_{k-1,Q}f\Vert^u_{L^u(Q)}\le cr_K^{n-d}\Vert f\Vert^u_{L^u{(25K\cap S,\mathcal{H}^d)}},$$
where  $K=Q(x_K,r_K)$. It remains to use 
\eqref{sset_rem}.
\end{proof}

Next we estimate local approximations of $\tilde{f}$
on cubes which are located  relatively far from the set $S$. Namely, assume
that a cube $K=Q(x_K,r_K)$ satisfies 
\begin{equation}\label{EqFarCubes}
0<{\rm diam}\,K\le {\rm dist}(x_K,S)/40\le 1024\Delta.
\end{equation}
Let $Q_K\in\mathcal{W}_S$ denote the Whitney cube which contains the center $x_K$ of $K$,
and let $a_K\in S$ be such
that $\|a_K-x_K\|_\infty=\mathrm{dist}(x_K,S)$. Define also $$\tilde{Q}_K:=Q(a_K,2{\rm dist}(x_K,S)),\qquad\hat{Q}_K:=25\tilde{Q}_K=Q(a_K,50{\rm dist}(x_K,S)).$$
By \eqref{EqFarCubes}, we have $K\subset\tilde{Q}_K$. Set
$$\mathcal{A}_K:=\{Q\in\mathcal{Q}_2(\tilde{Q}_K):\;{\rm diam}\,Q\le 4{\rm dist}(x_K,S)\}.$$

\begin{lem}\label{lem316}
For every cube $K$ satisfying (\ref{EqFarCubes}) and  every  $P_0\in \mathcal{P}_{k-1}$,
$$\mathcal{E}_k(\tilde{f},K)^u_{L^u(\R^n)}\le c_{n,u,k} \bigg(\frac{{\rm diam}\;K}{{\rm dist}(x_K,S)}\bigg)^{ku}|\tilde{Q}_K|^{-1}\sum_{Q\in\mathcal{A}_K}\Vert P_{k-1,Q}f-P_0\Vert^u_{L^u(Q)}.$$
\end{lem}

This lemma is proven in \cite[p. 1226]{Shvartsman}.

\begin{lem}\label{lem317}
Suppose that a cube $K$ satisfies both (\ref{EqFarCubes}) and ${\rm dist}(x_K,S)\le\Delta/4$. Then
$$\mathcal{E}_k(\tilde{f},K)_{L^u(\R^n)}\le c \bigg(\frac{{\rm diam}\;K}{{\rm dist}(x_K,S)}\bigg)^{k}\mathcal{E}_k(f,\hat{Q}_K)_{L^u(S)}.$$
Here the constant $c$ depends on $S,\Delta,n,u$ and $k$.
\end{lem}

\begin{proof}
Since ${\rm dist}(x_K,S)\le\Delta/4$, we have $\mathcal{A}_K\subset \mathcal{Q}_3(\tilde{Q}_K)$. Hence, by Lemma \ref{lem316}, for every $P_0\in \mathcal{P}_{k-1}$ we have
$$\mathcal{E}_k(\tilde{f},K)^u_{L^u(\R^n)}\le c \bigg(\frac{{\rm diam}\;K}{{\rm dist}(x_K,S)}\bigg)^{ku}|\tilde{Q}_K|^{-1}\sum_{Q\in\mathcal{Q}_3(\tilde{Q}_K)}\Vert P_{k-1,Q}f-P_0\Vert^u_{L^u(Q)}.$$
Notice that $\tilde{Q}_K$ is centered in $S$. Therefore,
by Lemma \ref{knine}, we have
\begin{align*}
&\sum_{Q\in\mathcal{Q}_3(\tilde{Q}_K)}\Vert P_{k-1,Q}f-P_0\Vert^u_{L^u(Q)}
\\&\le cr_{\tilde{Q}_K}^{n-d}\Vert f-P_0\Vert^u_{L^u(25\tilde{Q}_K\cap S,\mathcal{H}^d)}=cr_{\tilde{Q}_K}^{n-d}\Vert f-P_0\Vert^u_{L^u(\hat{Q}_K\cap S,\mathcal{H}^d)}.
\end{align*}
Let $P_0\in\mathcal{P}_{k-1}$ be a polynomial of the best approximation of $f$ on $\hat{Q}_K\cap S$, i.e.
\[
\Vert f-P_0\Vert^u_{L^u(\hat{Q}_K\cap S)} = \mathcal{H}^d(\hat Q_K\cap S)\mathcal{E}_k(f,\hat Q_K)_{L^u(S)}^u.
\]
Combining the estimates above, we find that
$$\mathcal{E}_k(\tilde{f},K)^u_{L^u(\R^n)}\le c \bigg(\frac{{\rm diam}\;K}{{\rm dist}(x_K,S)}\bigg)^{ku}\frac{\mathcal{H}^d(\hat{Q}_K\cap S)}{r_{\tilde{Q}_K}^{d}}\mathcal{E}_k(f,\hat{Q}_K)_{L^u(S)}^u.$$
Hence, the required estimate follows from \eqref{sset_rem}.
\end{proof}

\begin{lem}\label{lem318}
Suppose that a cube $K$ satisfies (\ref{EqFarCubes}), then
$$\mathcal{E}_k(\tilde{f},K)_{L^u(\R^n)}\le c \bigg(\frac{{\rm diam}\;K}{{\rm dist}(x_K,S)}\bigg)^{k}\bigg(\kint_{\hat{Q}_K\cap S}|f|^{u}\,d\mathcal{H}^d\bigg)^{\frac{1}{u}},$$
where the constant $c$ depends on $S,n,u,\Delta$ and $k$.
\end{lem}

\begin{proof}
Recall that if ${\rm diam}\,Q>\Delta$, then $P_{k-1,Q}f=0$, thus
$$
\sum_{Q\in\mathcal{A}_K}\Vert P_{k-1,Q}f\Vert^u_{L^u(Q)}\le
\sum_{Q\in\mathcal{Q}_2(\tilde{Q}_K)}\Vert P_{k-1,Q}f\Vert^u_{L^u(Q)}=\sum_{Q\in\mathcal{Q}_3(\tilde{Q}_K)}\Vert P_{k-1,Q}f\Vert^u_{L^u(Q)}.$$
By Lemma \ref{lem316} with $P_0=0$, we have
$$
\mathcal{E}_k(\tilde{f},K)^u_{L^u(\R^n)}\le c \bigg(\frac{{\rm diam}\;K}{{\rm dist}(x_K,S)}\bigg)^{ku}|\tilde{Q}_K|^{-1}\sum_{Q\in\mathcal{Q}_3(\tilde{Q}_K)}\Vert P_{k-1,Q}f\Vert^u_{L^u(Q)}.
$$
Hence, by Lemma \ref{knine} and \eqref{sset_rem},
we get
\begin{align*}
\mathcal{E}_k(\tilde{f},K)^u_{L^u(\R^n)} &\le c \bigg(\frac{{\rm diam}\;K}{{\rm dist}(x_K,S)}\bigg)^{ku}r_{\tilde{Q}_K}^{-d}\Vert f\Vert^u_{L^u(25\tilde{Q}_K\cap S)}\\&\le c \bigg(\frac{{\rm diam}\;K}{{\rm dist}(x_K,S)}\bigg)^{ku}\frac{1}{\mathcal{H}^d(\hat{Q}_K\cap S)}\Vert f\Vert^u_{L^u(\hat{Q}_K\cap S)}.
\end{align*}
\end{proof}

We are ready for the proof of Proposition \ref{th36}.
\begin{proof}
Let us first verify the identity \eqref{rest}.
Suppose that $f\in L^u_{\mathrm{loc}}(S)$. It suffices
to show that
\begin{equation}\label{vera}
\lim_{r\to 0+}\bigg|\kint_{Q(x,r)} \tilde f(y)dy - f(x)\bigg|=0\qquad \mathcal{H}^d\text{ a.e. }x\in S.
\end{equation}
Lemma \ref{Proposition310} and \eqref{sset_rem} imply
\begin{equation}\label{vas}
\bigg|\kint_{Q(x,r)} \tilde f(y)dy - f(x)\bigg|\le c\kint_{B(x,25\cdot\sqrt n\cdot r)\cap S} |f(y)-f(x)|d\mathcal{H}^d(y).
\end{equation}
By Remark 2.15 in \cite{Mat95},
the right hand side of \eqref{vas} tends to zero
for $\mathcal{H}^d$ almost every point $x\in S$ as $r\to 0$. 

Next we prove \eqref{norm}.
Let $K:=Q(x,t)$, where $x\in\R^n$ and $t>0$. Consider the following four cases {\bf I--IV}.

{\bf I}: $80t\le{\rm dist}(x,S)$ and $r^{(x,t)}\le\Delta$. 
By Lemma \ref{lem317}
$$\mathcal{E}_k(\tilde{f},K)_{L^u(\R^n)}\le c \bigg(\frac{{\rm diam}\,K}{{\rm dist}(x_K,S)}\bigg)^{k}\mathcal{E}_k(f,\hat{Q}_K)_{L^u(S)},$$
where $\hat{Q}_K=Q(a_x,50{\rm dist}(x,S))=K^{(x,t)}$. 
Since $t<80t\le {\rm dist}(x,S)$, we have $2{\rm dist}(x,S)^k\ge t^k+{\rm dist}(x,S)^k$, and in case {\bf I} the theorem is proven.

{\bf II}: $80t\le{\rm dist}(x,S)$ and $r^{(x,t)}\in (\Delta,1024\Delta]$. The reasoning in this case is the same as in case {\bf I}, except that we apply
Lemma \ref{lem318} instead of Lemma \ref{lem317}.

{\bf III}:
$80t>{\rm dist}(x,S)$ and $r^{(x,t)}\le\Delta$.  
Set $\hat{K}:=Q(a_x,81t)$. Since $K\subset\hat{K}$, we have
\begin{equation}\label{Case3}
\mathcal{E}_k(\tilde{f},K)_{L^u(\R^n)}\le 
\bigg(\frac{|\hat{K}|}{|K|}\bigg)^{1/u}\mathcal{E}_k(\tilde{f},\hat{K})_{L^u(\R^n)}=
c_{n,u} \mathcal{E}_k(\tilde{f},\hat{K})_{L^u(\R^n)}.
\end{equation}
Lemma \ref{th311} implies that
$\mathcal{E}_k(\tilde{f},\hat{K})_{L^u(\R^n)}\le c \mathcal{E}_k(f,25\hat{K})_{L^u(S)}$ 
and, by  \eqref{sset_rem}, we get
\begin{align*}
&\mathcal{E}_k(f,25\hat{K})_{L^u(S)}
 \le c_S\mathcal{E}_k(f,K^{(x,t)})_{L^u(S)}.
\end{align*}
It remains to note that $c_k t^k\ge t^k+{\rm dist}(x,S)^k$, and
the present case is proved.

{\bf IV}:
$80t>{\rm dist}(x,S)$ and $r^{(x,t)}\in (\Delta,1024\Delta]$. Reasoning as in the previous case we obtain (\ref{Case3}). By Lemma \ref{p312},
\[\mathcal{E}_k(\tilde{f},\hat{K})_{L^u(\R^n)}\le
\bigg(\kint_{\hat K}|\tilde f|^udx\bigg)^{1/u}\le
c\bigg(
\kint_{25\hat{K}\cap S}|f|^{u}\,d\mathcal{H}^d\bigg)^{\frac{1}{u}}.
\] 
Taking into account (\ref{Case3}) and the last inequality, we have
\begin{align*}
\mathcal{E}_k(\tilde{f},K)_{L^u(\R^n)}\le c\mathcal{E}_k(\tilde{f},\hat{K})_{L^u(\R^n)}
&\le c\bigg(\kint_{25\hat{K}\cap S}
|f|^{u}\,d\mathcal{H}^d\bigg)^{\frac{1}{u}}.
\end{align*}
Using \eqref{sset_rem}  and the fact that
$25\hat{K}\subset K^{(x,t)}$
finishes
the proof of case {\bf IV} and the theorem.
\end{proof}

\subsection{Global estimates for the extension}

Our next step is to establish Proposition \ref{fact}. 
It provides the key estimates for the proof of the extension theorems in Section \ref{ext_th_sec}.

Let $f\in L^u_{\mathrm{loc}}(S)$, $1\le u<\infty$ and  $k\in\N$. Suppose that $\tilde f=\mathrm{Ext}_{k,S}$ is defined as in \eqref{DefExtensionOperator} with
$
\Delta=16000.
$

Recall that
$\mathcal{D}$ denotes the family of closed dyadic cubes in $\R^n$ and
$\mathcal{D}_j$ stands for
the family of those dyadic cubes whose
side length is $2^{-j}$.  Denote also
\[
\mathcal{C}:=\{Q\in\mathcal{D}\,:\,\mathrm{dist}(x_Q,S)/320\le \ell(Q)\le 1\}.
\] 

For $i\in\Z$, we define
\begin{equation}\label{s_def}
S_i:=\{y\in\R^n\,:\,80\cdot 2^{-i} \le \mathrm{dist}(y,S)\le 16000\cdot 2^{-i}\}.
\end{equation}
Choose $m_0\in\N$ such that $2^{-m_0}< 1/200$. Then
\begin{equation}\label{s_disj}
S_i\cap S_{i'}=\emptyset,
\end{equation}
if $i-i'\equiv 0\,\mathrm{mod}\,m_0$ and $i\not=i'$. 

Let
$\mathcal{F}_i$, $i\in\Z$,  
be a finite or countable family of cubes
of the form $Q(x_Q,2^{-i}\Delta)$, $x_Q\in S$,
with the following properties:
\begin{equation}\label{sov}
\begin{split}
&S\subset \bigcup_{Q\in \mathcal{F}_i} Q;\\
&\sup_{Q\in\mathcal{F}_i}\,\, {\rm card}\,\{R\in \mathcal{F}_i\,:\,2R\cap 2Q\not=\emptyset\}\le q_0,
\end{split}
\end{equation}
where $q_0 \in\N$ is a constant.
The existence of  $\mathcal{F}_i$ follows
from the $5r$-covering theorem, see e.g. \cite[p. 23]{Mat95}. The lemma below is a direct consequence of \eqref{sov}.
\begin{lem} 
Let $i\in\Z$. Then
we have the following decomposition
\begin{equation}\label{fkaks}
2\mathcal{F}_i:=\{2Q\,:\,Q\in \mathcal{F}_i\}=\bigcup_{q=1}^{q_0} 2\mathcal{F}_i^q,
\end{equation}
where each $2\mathcal{F}_i^q$ is a disjoint family (possibly empty). 
\end{lem}

Using the introduced notation we can formulate the main result of this section.
\begin{prop}\label{fact}
Let $S$ be a $d$-set in $\R^n$ with $n-1<d<n$.
Let $k\in\N$, $1\le u <\infty$ and $f\in L^u_{\mathrm{loc}}(S)$. If $j\in \N_0$, then
\begin{equation}\label{cubetrans}
\mathcal{E}_k(\tilde f,Q(x,2^{-j}))_{L^u(\R^n)}
\le c\sum_{Q\in\mathcal{D}_j} \chi_Q(x) \mathcal{E}_k(\tilde f,4Q)_{L^u(\R^n)},\quad x\in\R^n.
\end{equation}
Furthermore, we have 
\begin{equation}\label{eest}
\sum_{\substack{Q\in \mathcal{D}_j\cap \mathcal{C}}}\chi_Q \mathcal{E}_k(\tilde f,4Q)_{L^u(\R^n)}\le 
c\sum_{K\in\mathcal{F}_j}\chi_{2K}\mathcal{E}_{k}(f,2K)_{L^u(S)};
\end{equation}
and
\begin{equation}\label{test}
\sum_{\substack{Q\in \mathcal{D}_j\setminus \mathcal{C}}}\chi_Q \mathcal{E}_k(\tilde f,4Q)_{L^u(\R^n)}\le 
c\sum_{i=-10}^{j-1} 2^{k(i-j)}\sum_{K\in\mathcal{F}_i}\chi_{2K\cap S_i}\mathcal{E}_{k(i)}(f,2K)_{L^u(S)}.
\end{equation}
Here we denoted $k(i)=k$ if $i\ge 0$ and $k(i)=0$ otherwise, and the constant
$c$ depends on $S,n,u$ and $k$.
\end{prop}

\begin{proof}
The ideas of the proof are similar to the ones used in Theorem~5.3 in \cite{Shvartsman}. 

First, note that \eqref{cubetrans} follows from the monotonicity property 
\eqref{eqMonotonyOfLocalApproxSset}
of local approximation. 

Next we partition $\mathcal{D}_j$
into disjoint families $\mathcal{D}_j^i$, $i\le j$, 
in the following way: define $\mathcal{D}_j^j= \mathcal{D}_j\cap \mathcal{C}$ 
and, for each integer $i<j$, 
\begin{equation}\label{i_val}
\mathcal{D}_j^i=\big\{Q\in\mathcal{D}_j\setminus \mathcal{C}:\;320\cdot 2^{-i-1}
<{\rm dist}(x_Q,S)\le 320\cdot 2^{-i}\big\}.
\end{equation}
Then
\begin{equation}\label{ssaf}
\mathcal{D}_j\setminus \mathcal{C} =\{Q\in \mathcal{D}_j\setminus \mathcal{C}\,:\,\mathrm{dist}(x_Q,S)> 320\cdot 2^{-j}\}=\bigcup_{i<j} \mathcal{D}_j^i
\end{equation}
and
$\mathcal{D}_j=\cup_{i\le j}\mathcal{D}_j^i$.

Recall that for $4Q=Q(x_Q,2^{-j+1})$, where $Q\in\mathcal{D}_j$, we denote
\begin{equation}\label{beg}
\begin{split}
r^{(x_Q,2^{-j+1})}&=50\max(80\cdot 2^{-j+1},{\rm dist}(x_Q,S));\\K^{(x_Q,2^{-j+1})}&=Q(a_{x_Q},r^{(x_Q,2^{-j+1})}),
\end{split}
\end{equation}
where $a_{x_Q}$ is a point on $S$ nearest to $x_Q$.

Let $i<j$ and $Q\in \mathcal{D}_j^i$. 
Since $\mathcal{F}_i$ is a cover of $S$, we can choose a cube $K_Q\in \mathcal{F}_i$ for
which $a_{x_Q}\in K_Q$.
Observe that 
\begin{equation}\label{q_incl}
Q\subset Q(a_{x_Q},2^{-j-1}+{\rm dist}(x_Q,S))\subset Q(a_{x_Q},2^{-i}\Delta/25)\subset 2K_Q.
\end{equation}
Furthermore, by \eqref{ssaf}, we have
$80\cdot 2^{-j+1} < 320 \cdot 2^{-j} <\mathrm{dist}(x_Q,S)$ and 
\begin{equation}\label{isop}
r^{(x_Q,2^{-j+1})}=50{\rm dist}(x_Q,S)
\in 
\begin{cases}
(0,\Delta],\qquad &0\le i<j;\\
(\Delta,1024\Delta],\qquad &-10\le i<0.
\end{cases}
\end{equation}

Let us show that for every $Q\in \mathcal{D}_j^i$, $i< j$,
\begin{equation}\label{eqModulusOfCont}
\mathcal{E}_k(\tilde{f},4Q)_{L^u(\R^n)}\le 
c2^{k(i-j)}
\begin{cases}
\mathcal{E}_k(f,2K_Q)_{L^u(S)},\quad &0\le i< j;\\
\mathcal{E}_0(f,2K_Q)_{L^u(S)},\quad &i<0.
\end{cases}
\end{equation}
First recall that, by \eqref{supp}, 
$\tilde f(x)=0$ if $\mathrm{dist}(x,S)>8\Delta=125\cdot 2^{10}$.
On the other hand,
$$\mathrm{dist}(x,S)\ge 80\cdot 2^{-i},\qquad x\in 4Q.$$
Hence, 
\begin{equation}\label{whatis}
\mathcal{E}_k(\tilde f,4Q)_{L^u(\R^n)}=0,\qquad i< -10.
\end{equation}
If $i\ge -10$, by \eqref{isop} and Proposition \ref{th36}, we have
$$
\mathcal{E}_k(\tilde{f},4Q)_{L^u(\R^n)}\le c \frac{2^{-jk}}{{\rm dist}(x_Q,S)^k}\begin{cases}
\mathcal{E}_k(f,K^{(x_Q,2^{-j+1})})_{L^u(S)},\quad &0\le i< j;\\
\mathcal{E}_0(f,K^{(x_Q,2^{-j+1})})_{L^u(S)},\quad &-10\le i<0.
\end{cases}
$$
Here $K^{(x_Q,2^{-j+1})}=Q(a_{x_Q},50{\rm dist}(x_Q,S))$.  By \eqref{i_val},
$$K^{(x_Q,2^{-j+1})}\subset Q(a_{x_Q},2^{-i}\Delta)\subset 2K_Q\text{ and } \ell(K_Q)/4< 50{\rm dist}(x_Q,S)\le \ell(K_Q)/2.$$
Thus, \eqref{eqModulusOfCont} follows 
from \eqref{eqMonotonyOfLocalApproxSset} and the estimates
above.

By \eqref{ssaf},  \eqref{eqModulusOfCont} and \eqref{whatis}, we have
\begin{equation}\label{refids}
\begin{split}
&\sum_{Q\in \mathcal{D}_j\setminus \mathcal{C}} \chi_Q\mathcal{E}_k(\tilde{f},4Q)_{L^u(\R^n)}
\le c\sum_{i=-10}^{j-1} 
2^{k(i-j)}\sum_{Q\in\mathcal{D}_j^i}\chi_Q \mathcal{E}_{k(i)}(f,2K_Q)_{L^u(S)}.
\end{split}
\end{equation}
Now we will rewrite the right hand side of \eqref{refids} in terms of
the cubes from $\mathcal{F}_i$.
To this end, first notice that
\begin{align*}
\sum_{Q\in\mathcal{D}_j^i}\chi_Q \mathcal{E}_{k(i)}(f,2K_Q)_{L^u(S)}=
\sum_{K\in\mathcal{F}_i}\sum_{\substack{Q\in\mathcal{D}_j^i\\K_Q=K}}\chi_Q\mathcal{E}_{k(i)}(f,2K)_{L^u(S)}.
\end{align*}
On the other hand, if $Q\in \mathcal{D}_j^i$, $i<j$, then $Q\subset S_i$ and,
by \eqref{q_incl}, we have 
$Q\subset 2K_Q\cap S_i$.
Note also that every point in $\R^n$ is covered  by at most $c_n$ cubes from $\mathcal{D}^i_j$. Hence,
\begin{equation}\label{EqOmeqa_I>0}
\sum_{Q\in\mathcal{D}_j^i}\chi_Q \mathcal{E}_{k(i)}(f,2K_Q)_{L^u(S)}\le  c\sum_{K\in\mathcal{F}_i}\chi_{2K\cap S_i}\mathcal{E}_{k(i)}(f,2K)_{L^u(S)},\quad i< j.
\end{equation}
Combining \eqref{refids} and \eqref{EqOmeqa_I>0} finishes the proof of  estimate  \eqref{test}.

Estimate
\eqref{eest} is obtained by similar reasoning, omitting the arguments about the sets $S_i$.
\end{proof}

\section{Extension theorems}\label{ext_th_sec}

In this section we prove the extension theorems for Besov
and Triebel--Lizorkin spaces. The proofs of these extension results share common technical estimates which have already been obtained, see Proposition 
\ref{fact}.

\subsection{Besov spaces}
The extension theorem for Besov spaces is the following.
\begin{thm}\label{bes_est}
Let $S$ be a $d$-set in $\R^n$ with $n-1<d<n$.
Suppose that $k\in\N$, $0<\alpha<k$, $1\le p,q<\infty$ and $1\le u\le p$.
Then, for every function $f\in L^p(S)$,
\begin{equation}\label{ExtTheoremBesovSpaces}
\begin{split}
&\bigg(\sum_{j=0}^\infty \big(2^{\alpha j}\|\mathcal{E}_k(\tilde f,Q(\cdot,2^{-j}))_{L^u(\R^n)}\|_{L^p(\R^n)}\big)^q\bigg)^{1/q} + \|\tilde f\|_{L^p(\R^n)}\\
&\le c\Bigg\{\bigg(\sum_{j=0}^\infty \big(
2^{j(\alpha-(n-d)/p)}\|\mathcal{E}_k(f,Q(\cdot,2^{-j}))_{L^u(S)}\|_{L^p(S)}\big)^q\bigg)^{1/q}+\Vert f\Vert_{L^p(S)}\Bigg\}.
\end{split}
\end{equation}
Here the constant $c$ depends on $S,n,u,k,p,q$ and $\alpha$.
\end{thm}

Before the proof of the theorem we consider the following proposition.
\begin{prop}\label{th53}
Let $S$ be a $d$-set in $\R^n$ with $n-1<d<n$.
Let $1\le u\le p<\infty$ and $f\in L^p(S)$. Then, for every $j\in\N_0$, 
\begin{equation}\label{six2}
\begin{split}
&\Vert\mathcal{E}_k(\tilde{f},Q(\cdot,2^{-j}))_{L^u(\R^n)}\Vert_{L^p(\mathbb{R}^n)}\\
&\le c2^{-jk}\bigg\{\bigg(\sum_{i=0}^j\big(2^{i(k-(n-d)/p)}\Vert\mathcal{E}_k(f,Q(\cdot,2^{-i}))_{L^u(S)}\Vert_{L^p(S)}\big)^p\bigg)^\frac{1}{p}+\Vert f\Vert_{L^p(S)}\bigg\},
\end{split}
\end{equation}
where the constant $c$ depends on $S,n,u,k$ and $p$.
\end{prop}

\begin{proof}
By Proposition \ref{fact},
\[
\Vert\mathcal{E}_k(\tilde{f},Q(\cdot,2^{-j}))_{L^u(\R^n)}\Vert_{L^p(\mathbb{R}^n)}\le I+II+III,
\]
where
\begin{align*}
&I:=
\bigg\|\sum_{i=0}^{j-1} 2^{k(i-j)}\sum_{K\in\mathcal{F}_i}\chi_{2K\cap S_i}\mathcal{E}_{k}(f,2K)_{L^u(S)}\bigg\|_p;\\
&II:=\bigg\|\sum_{i=-10}^{-1} 2^{k(i-j)}\sum_{K\in\mathcal{F}_i}\chi_{2K\cap S_i}\mathcal{E}_{0}(f,2K)_{L^u(S)}\bigg\|_p;\\
&III:=\bigg\| \sum_{K\in\mathcal{F}_j}\chi_{2K}\mathcal{E}_{k}(f,2K)_{L^u(S)}\bigg\|_p.
\end{align*}
Let us first focus on the term $I$. By \eqref{fkaks}, we write
\begin{align*}
I^p&=2^{-jkp}\int_{\R^n} \bigg|
\sum_{m=0}^{m_0-1}\sum_{\substack{i=0\ldots j-1\\ i\equiv m\,\mathrm{mod}\,m_0}} 2^{ik}
\sum_{q=1}^{q_0}
\sum_{Q\in 2\mathcal{F}_i^q}\chi_{Q\cap S_i}(x)\mathcal{E}_{k}(f,Q)_{L^u(S)}\bigg|^pdx\\
&\le cm_0^pq_0^p 2^{-jkp}\sum_{m=0}^{m_0-1}  \sum_{q=1}^{q_0}  
\underbrace{\int_{\R^n} \bigg|
\sum_{\substack{i=0\ldots j-1\\i\equiv m\,\mathrm{mod}\,m_0}}  2^{ik}\sum_{Q\in 2\mathcal{F}_i^q}\chi_{Q\cap S_i}(x)\mathcal{E}_{k}(f,Q)_{L^u(S)}\bigg|^pdx}_{=:I_{m,q}^p}.
\end{align*}
Fix $m\in \{0,\ldots,m_0-1\}$ and $q\in\{1,\ldots,q_0\}$.
Observe that
 the family 
\[\{Q\cap S_i\,:\,i\ge 0,\,i\equiv m\,\mathrm{mod}\,m_0,\,Q\in 2\mathcal{F}_i^q\}\] is disjoint. This fact follows from \eqref{s_disj} and the disjointness 
of the family $2\mathcal{F}_j^q$. Hence, we have
\begin{align*}
I_{m,q}^p
&= \int_{\R^n}
\sum_{\substack{i=0\ldots j-1\\i\equiv m\,\mathrm{mod}\,m_0}} 2^{ikp}\sum_{Q\in 2\mathcal{F}_i^q}\chi_{Q\cap S_i}(x)\mathcal{E}_{k}(f,Q)_{L^u(S)}^pdx
\\&\le \sum_{i=0}^{j-1} 2^{ikp}\sum_{Q\in 2\mathcal{F}_i^q}{|Q|}\mathcal{E}_{k}(f,Q)_{L^u(S)}^p.
\end{align*}
For every $Q\in 2\mathcal{F}_i$, $i\ge 0$, we have $\ell(Q)= 2^{-i+2}\Delta\le 4 \Delta$ and, by \eqref{sset_rem},
$$\mathcal{H}^d(Q\cap S)\geq c (2^{-i}\Delta)^d.$$ 
Hence,
\begin{equation}\label{bo}
\sum_{Q\in 2\mathcal{F}_i^q}{|Q|}\mathcal{E}_{k}(f,Q)_{L^u(S)}^p \le c2^{i(d-n)}\sum_{Q\in 2\mathcal{F}_i^q}\mathcal{H}^d(Q\cap S)\mathcal{E}_{k}(f,Q)_{L^u(S)}^p.
\end{equation}
Let $i_\Delta\in \N$ such that $2^{i_\Delta-3}< \Delta\le 2^{i_\Delta-2}$. Note that $Q\subset Q(x,2^{-i+i_\Delta})$ for every $x\in Q$. Then, by 
\eqref{eqMonotonyOfLocalApproxSset} and the fact that $2\mathcal{F}_i^q$ is a disjoint family, we have
\begin{equation}\label{i_iso}
\sum_{Q\in 2\mathcal{F}_i^q}\mathcal{H}^d(Q\cap S)\mathcal{E}_{k}(f,Q)_{L^u(S)}^p
\le c\|\mathcal{E}_k(f,Q(\cdot,2^{-i+i_\Delta}))_{L^u(S)}\|_{L^p(S)}^p.
\end{equation}
Notice also that, by Fubini's theorem and \eqref{sset_rem},
\[\|\mathcal{E}_k(f,Q(\cdot,2^{-i+i_\Delta}))_{L^u(S)}\|_{L^p(S)}^p\le c\|f\|_{L^p(S)}^p.\]
Using the last estimate 
for the indices $i=0,\ldots,i_\Delta-1$, and \eqref{i_iso} for indices $i\ge i_\Delta$, we obtain
\[
I_{m,q}^p \le c\bigg\{\sum_{i=0}^{j} \big(
2^{i(k-(n-d)/p)}\|\mathcal{E}_k(f,Q(\cdot,2^{-i}))_{L^u(S)}\|_{L^p(S)}\big)^p+\Vert f\Vert_{L^p(S)}^p\Bigg\}.
\] 

Performing similar computations for the term $III$, we see that $I+III$ is bounded (up to a constant) by the right hand side of
\eqref{six2}.

Let us estimate the remaining term $II$.
By the triangle-inequality, we get
\begin{align*}
II&\le c2^{-jk}\sum_{i=-10}^{-1}\sum_{q=1}^{q_0}\bigg( \int_{\R^n} \bigg(2^{ik}\sum_{Q\in 2\mathcal{F}_i^q}\chi_{Q\cap S_i}(x)\mathcal{E}_0(f,Q)_{L^u(S)}\bigg)^p  dx\bigg)^{1/p}.
\end{align*}
Recall that the family $2\mathcal{F}_i^q$ is disjoint. Hence,
\begin{align*}
II\le &c2^{-jk}\sum_{i=-10}^{-1}\sum_{q=1}^{q_0} \bigg(2^{ikp}\sum_{Q\in 2\mathcal{F}_i^q}|Q|\mathcal{E}_0(f,Q)_{L^u(S)}^p\bigg)^{1/p}
\end{align*}
Let $i\in \{-10,\ldots,-1\}$ and $q\in \{1,\ldots,q_0\}$.
For every $Q\in 2\mathcal{F}_i^q$ we have $\ell(Q)= 2^{-i+2}\Delta\le 4096 \Delta$ and, by \eqref{sset_rem},
$$|Q|\le c 2^{-i(n-d)}\mathcal{H}^d(Q\cap S)\le  c\mathcal{H}^d(Q\cap S).$$
Then by H\"{o}lder's inequality, we get
\begin{align*}
|Q|\mathcal{E}_0(f,Q)^p_{L^u(S)}
&=|Q|\bigg(\frac{1}{\mathcal{H}^d(Q\cap S)}\int_{Q\cap S}|f|^u\,d\mathcal{H}^d\bigg)^{p/u}
\le c\int_{Q\cap S}|f|^{p}\,d\mathcal{H}^d.
\end{align*}
Since $2\mathcal{F}_i^q$ is disjoint, we see that
\begin{align*}
II^p\le
c2^{-jkp}\sum_{Q\in 2\mathcal{F}_i^q}\;\int_{Q\cap S}|f|^{p}\,d\mathcal{H}^d\le  c2^{-jkp}\|f\|_{L^p(S)}^p.
\end{align*}
\end{proof}

We are now ready for the proof of Theorem \ref{bes_est}.

\begin{proof}
Let us first estimate $\|\tilde f\|_{L^p(\R^n)}$. Recall that, by \eqref{supp},
$
\mathrm{supp}\,\tilde f\subset \bigcup_{x\in S} Q(x,8\Delta)$. 
By $5r$-covering theorem there are points $x_m\in S$, $m=1,2,\ldots$ (possibly a finite number of them) such that the cubes $Q(x_m,8\Delta)$ are pairwise disjoint and
\[\mathrm{supp}\,\tilde f\subset \bigcup_{m=1}^\infty Q(x_m,40 \Delta).\]
By Lemma \ref{p312}, we obtain
\begin{align*}
\|\tilde f\|_{L^p(\R^n)}^p &\le \sum_{m=1}^\infty \int_{Q(x_m,40\Delta)} |\tilde f(x)|^pdx\\
&\le c\sum_{m=1}^\infty \frac{|Q(x_m,40 \Delta)|}{\mathcal{H}^d(Q(x_m,1000\Delta)\cap S)}\int_{Q(x_m,1000\Delta)\cap S} |f(x)|^pd\mathcal{H}^d(x)\\&\le c\|f\|_{L^p(S)}^p,
\end{align*}
where the last inequality follows from 
\eqref{sset_rem} and the fact that \[\sum_{m=1}^\infty \chi_{Q(x_m,1000 \Delta)}\le c_n.\]

Next we focus on the remaining series at the left hand side of 
\eqref{ExtTheoremBesovSpaces}.
For $j\ge 0$, denote
\begin{equation*}
A_j :=\sum_{i=0}^j a_i := \sum_{i=0}^j\big(2^{i(k-(n-d)/p)}\Vert\mathcal{E}_k(f,Q(\cdot,2^{-i}))_{L^u(S)}\Vert_{L^p(S)}\big)^p.
\end{equation*}
Then by Proposition \ref{th53}, we have
\begin{align*}
&\bigg(\sum_{j=0}^\infty \big(2^{j\alpha}\|\mathcal{E}_k(\tilde f,Q(\cdot,2^{-j}))_{L^u(\R^n)}\|_{L^p(\R^n)}\big)^q\bigg)^{1/q}
\\&\le 
c\bigg(\sum_{j=0}^\infty \big(2^{j(\alpha-k)}
\big\{A_j^{1/p} + \|f\|_{L^p(S)}\big\}\big)^q\bigg)^{1/q}
\\&\le c\bigg(\sum_{j=0}^\infty 2^{jq(\alpha-k)}
A_j^{q/p}\bigg)^{1/q} + c\|f\|_{L^p(S)}.
\end{align*}
For the last inequality we use the triangle inequality and the assumption that $q(\alpha-k)<0$.

By \eqref{leindler},
\begin{equation*}
\begin{split}
\bigg(\sum_{j=0}^\infty 2^{jq(\alpha-k)}
A_j^{q/p}\bigg)^{1/q}&=\bigg(\sum_{j=0}^\infty 2^{jq(\alpha-k)}
\bigg(\sum_{i=0}^j a_i\bigg)^{q/p}\bigg)^{1/q}\\&\le 
c\bigg(\sum_{j=0}^\infty 2^{jq(\alpha-k)}
 a_j^{q/p}\bigg)^{1/q}\\
&=c\bigg(\sum_{j=0}^\infty 
2^{jq(\alpha-(n-d)/p)}\Vert\mathcal{E}_k(f,Q(\cdot,2^{-j}))_{L^u(S)}\Vert_{L^p(S)}^q\bigg)^{1/q}.
\end{split}
\end{equation*}
\end{proof}

\subsection{Triebel--Lizorkin spaces}

For Triebel--Lizorkin spaces we have the following  extension theorem.
\begin{thm}\label{thf}
Let $S$ be a $d$-set in $\R^n$ with $n-1<d<n$.
 Let  $1< p<\infty$, $k\in\N$ and $0<\alpha<k$.
Then, for every function $f\in L^p(S)$, we have
\begin{align*}
&\bigg\| \sum_{j=0}^\infty 2^{j\alpha} \mathcal{E}_k(\tilde f,Q(\cdot,2^{-j}))_{L^1(\R^n)} \bigg\|_p+\|\tilde f\|_{L^p(\R^n)}\\
&\le c\Bigg\{\bigg(\sum_{j=0}^\infty \big(
2^{j(\alpha-(n-d)/p)}\|\mathcal{E}_k(f,Q(\cdot,2^{-j}))_{L^1(S)}\|_{L^p(S)}\big)^p\bigg)^{1/p}+\Vert f\Vert_{L^p(S)}\Bigg\}.
\end{align*}
Here the constant $c$ depends on $S,n,u,k,p,q$ and $\alpha$.
\end{thm}

\begin{proof}
The same reasoning as in the proof of Theorem \ref{bes_est} gives
$\|\tilde f\|_{L^p(\R^n)}\le c\|f\|_{L^p(S)}$.
By \eqref{cubetrans} in Proposition \ref{fact}, it is enough  to estimate the quantities
\begin{align*}
\bigg\| \sum_{\substack{Q\in\mathcal{C}}} \chi_Q r_Q^{-\alpha} \mathcal{E}_k(\tilde f,4Q)_{L^1(\R^n)}\bigg\|_p
+ \bigg\| \sum_{j=0}^\infty\sum_{\substack{Q\in\mathcal{D}_j\setminus \mathcal{C}}} \chi_Q r_Q^{-\alpha} \mathcal{E}_k(\tilde f,4Q)_{L^1(\R^n)}\bigg\|_p=:A+B.
\end{align*}
In order to estimate the first
term we take into account Remark \ref{porous_lem} and use 
Theorem \ref{reverse}. Note also that cubes from $\mathcal{D}_j$ have mutually disjoint interiors. Thus, we have
\begin{equation}\label{erot}
\begin{split}
A&\le
\bigg\| \bigg(\sum_{\substack{Q\in\mathcal{C}}} (\chi_Q r_Q^{-\alpha} \mathcal{E}_k(\tilde f,4Q)_{L^1(\R^n)})^p\bigg)^{1/p}\bigg\|_p\\
&\le c\bigg(\sum_{j=0}^\infty \int_{\R^n} \bigg(\sum_{Q\in\mathcal{D}_j} \chi_Q(x) r_Q^{-\alpha} \mathcal{E}_k(\tilde f,4Q)_{L^1(\R^n)}\bigg)^pdx\bigg)^{1/p}.
\end{split}
\end{equation}
Let $x\in Q\in\mathcal{D}_j$. Since $4Q\subset Q(x,8r_Q)$, by 
\eqref{eqMonotonyOfLocalApproxSset},
\begin{align*}
&\mathcal{E}_k(\tilde f,4Q)_{L^1(\R^n)}\le c\mathcal{E}_k(\tilde f,Q(x,8r_Q))_{L^1(\R^n)}= c\mathcal{E}_k(\tilde f,Q(x,2^{-j+2}))_{L^1(\R^n)}.
\end{align*}
As a consequence, we have
\[
\sum_{Q\in\mathcal{D}_j} \chi_Q(x) r_Q^{-\alpha} \mathcal{E}_k(\tilde f,4Q)_{L^1(\R^n)}
\le c2^{(j-2)\alpha}\mathcal{E}_k(\tilde f,Q(x,2^{-j+2}))_{L^1(\R^n)},\qquad x\in\R^n.
\]
Combining  \eqref{erot} and the last inequality, we have 
\[
A\le c\bigg\{\bigg(\sum_{j=0}^\infty \big(2^{j\alpha}\|\mathcal{E}_k(\tilde f,Q(\cdot,2^{-j}))_{L^1(\R^n)}\|_p\big)^p\bigg)^{1/p}+\| \tilde f\|_p\bigg\}.
\] 
Here for the indices $j=0$ and $j=1$ we use the estimate \[\|\mathcal{E}_k(\tilde f,Q(\cdot ,2^{-j+2}))_{L^1(\R^n)}\|_p^p\le 
c\| \tilde f\|_p^p.\] 
Applying Theorem \ref{bes_est} gives the required upper bound for term $A$.

To estimate term $B$, we use Proposition \ref{fact},
\begin{align*}
B&\le c\bigg\|   \sum_{j=0}^\infty 2^{j(\alpha-k)}\sum_{i=0}^{j} 2^{ik}\sum_{K\in\mathcal{F}_i}\chi_{2K\cap S_i}\mathcal{E}_{k}(f,2K)_{L^1(S)}   \bigg\|_p\\
&\qquad +c\bigg\|\sum_{j=0}^\infty 2^{j(\alpha-k)}\sum_{i=-10}^{-1} 2^{ik}\sum_{K\in\mathcal{F}_i}\chi_{2K\cap S_i}\mathcal{E}_{0}(f,2K)_{L^1(S)}\bigg\|_p=:B_1+B_2.
\end{align*}
By \eqref{leindler}, we get
\begin{align*}
&\sum_{j=0}^\infty 2^{j(\alpha-k)}\sum_{i=0}^{j} 2^{ik}\sum_{K\in\mathcal{F}_i}\chi_{2K\cap S_i}\mathcal{E}_{k}(f,2K)_{L^1(S)}   \\&\le c\sum_{j=0}^\infty 2^{j\alpha}\sum_{K\in\mathcal{F}_j}\chi_{2K\cap S_j}\mathcal{E}_{k}(f,2K)_{L^1(S)}.
\end{align*}
The arguments used for estimating term $I$ in the proof of Proposition \ref{th53} show that
\[
B_1^p\le c\bigg\{\sum_{j=0}^\infty \big(
2^{j(\alpha-(n-d)/p)}\|\mathcal{E}_k(f,Q(\cdot,2^{-j}))_{L^1(S)}\|_{L^p(S)}\big)^p+\Vert f\Vert_{L^p(S)}^p\Bigg\}.
\] 

Let us now consider term $B_2$. Since $\alpha-k<0$, we clearly  have 
\[
B_2\le c\bigg\|\sum_{i=-10}^{-1} 2^{ik}\sum_{K\in\mathcal{F}_i}\chi_{2K\cap S_i}\mathcal{E}_{0}(f,2K)_{L^1(S)}\bigg\|_p.
\]
Repeating the estimate for term $II$ in the proof of Proposition \ref{th53}, we see that $B_2\le c\|f\|_{L^p(S)}$.
\end{proof}

\section{Remark on the norm equivalence in Besov spaces}
Let us now complete 
the proof of Theorem \ref{beseq}. Recall that the norms in Besov space $B^\alpha_{pq}(S)$ defined by \eqref{DefinitionBesovNormIntegral}
involves two additional parameters $k$ and $u$. In Section~\ref{besov_sec} we showed that the definition does not depend on the choice of 
$k\in\N$ when $k>\alpha$, see Proposition \ref{noreq}.  

The following proposition shows that the definition of
Besov spaces on $d$-sets is independent
of the parameter $u$, as long as $1\le u\le p$. 

\begin{prop}\label{equirem}
Let $S$ be a $d$-set and $n-1<d<n$. Suppose that
$1\le p,q<\infty$ and $k>\alpha + (n-d)/p$ is an integer.
Then
Definition \ref{besov} of Besov space
$B_{pq}^{\alpha}(S)$ does not depend on the parameter
$1\le u\le p$ and the norms
corresponding to any pair $1\le u,u'\le p$ are
equivalent.
\end{prop}

\begin{proof}
Fix $k\in\N$ such that $k>\alpha+(n-d)/p$.
By Remark \ref{diskr}, we can replace integrals by sums in the definition of Besov spaces. Hence, if
we denote
\[
N_u(f):=\bigg(\sum_{j=0}^\infty \big(
2^{j\alpha}\|\mathcal{E}_k(f,Q(\cdot,2^{-j}))_{L^u(S)}\|_{L^p(S)}\big)^q\bigg)^{1/q}+\Vert f\Vert_{L^p(S)}
\]
for $f\in B_{pq}^{\alpha}(S)$, it suffices to verify that 
\begin{equation}\label{EquivalenceForU}
N_p(f)\le cN_1(f),
\end{equation}where
$c>0$ is independent of $f$. 
Indeed, the inequality 
$N_u(f)\le cN_{u'}(f)$ for any $1\le u,u'\le p$ simply follows from \eqref{EquivalenceForU}
and H\"older's inequality.

Denote $\tilde f:=\mathrm{Ext}_{k,S}(f)$ and recall that, by Proposition \ref{th36}, $\tilde f|_S=f$.
By Theorem \ref{restrB} we obtain the estimate \eqref{RestrictionForBesovSpace}.
To make a transition from $u=p$ to $u=1$ in the right hand side of \eqref{RestrictionForBesovSpace} 
we use an equivalence result for the norms of classical Besov spaces on $\R^n$, see Theorem 1.7.3 in \cite{Triebel2}.
This leads to the estimate
\begin{align*}
&N_p(f)\\&\le c\bigg\{\bigg(\sum_{j\ge 0} \big(2^{j (\alpha+(n-d)/p)}\Vert\mathcal{E}_k(\tilde f,Q(\cdot,2^{-j}))_{L^1(\R^n)}\Vert_{L^p(\R^n)}\big)^q\bigg)^{1/q}+\|\tilde f\|_{L^p(\R^n)}\bigg\},
\end{align*}
where the constant $c>0$ is independent of $f$.
Finally, by Theorem \ref{bes_est}, we see that the
right hand side of the last inequality is bounded (up to constant) by $N_1(f)$.
\end{proof}

By combining
Proposition
\ref{noreq} and Proposition \ref{equirem}, we obtain a proof
for Theorem~\ref{beseq}, which states that  Besov spaces on $d$-sets are independent
of the parameters $u$ and $k$, as long as $1\le u\le p$ and $k>\alpha$ is an integer.

\section*{Acknowledgements}
The authors are grateful to J. Kinnunen for 
having pointed out the problem, and for his valuable comments.

L. Ihnatsyeva was supported by the Academy of Finland.
A.V. V\"ah\"akangas was supported by the Academy of Finland CoE in Analysis and Dynamics research.

\bibliographystyle{elsarticle-num}

\medskip

\noindent Addresses:

\medskip
\noindent L.I.: Department of Mathematics, P.O. Box 11100, 
FI-00076 Aalto University, Finland. \\
\noindent 
e-mail: {\tt lizaveta.ihnatsyeva@aalto.fi}

\medskip
\noindent A.V.V.: Department of Mathematics and Statistics,
P.O. Box 68, FI-00014 University of Helsinki, Finland. \\
\noindent 
e-mail: {\tt antti.vahakangas@helsinki.fi}\\

\end{document}